\documentclass[12pt]{amsart}
\usepackage{ShaharStyle}
\title{Monodromic Hecke Modules for Symmetric Spaces}
\author{Shahar Dagan}
\thanks{Email: \texttt{shahar.dagan@weizmann.ac.il}}
\begin{document}
\begin{abstract}
    Let $G$ be a connected reductive group over the algebraic closure
    of a finite field of odd characteristic, with smooth connected center.
    Let $F$ be a Frobenius endomorphism and $\theta$ an involution
    commuting with $F$. Put $H\coloneqq G^\theta$ and $H^\circ$ its identity component.
    We compute the action of the simple Hecke generators on the
    standard basis of the $H^\circ$-equivariant monodromic Hecke module
    as torus monodromy varies over a Weyl-group orbit.
    The formulas include generators that change monodromy and recover
    the Lusztig--Vogan formulas when monodromy is trivial.
    Specialization at $u^{1/2}=-1$ yields a Weyl-group representation
    with a filtration by anti-invariant torus rank.
    The associated graded carries a signed cross action.
\end{abstract}
\maketitle
\section{Introduction and main results}
\label{sec:introduction}
\setcounter{tocdepth}{2}
\tableofcontents

\subsection{Statement of main results}
\label{subsec:Statement-of-main-results}

Let $G$ be a connected reductive group over
$k\coloneqq\overline{\bF}_q$, where $q$ is odd.
Assume that the scheme-theoretic center $Z(G)$ is smooth and connected.
Fix a Frobenius endomorphism $F$ and an involution $\theta$ commuting
with $F$. Put $H\coloneqq G^\theta$, and write $H^\circ$ for its
identity component.
Fix a prime $\ell\ne\operatorname{char}(k)$.
All local systems have coefficients in $\overline{\bQ}_\ell$.

Choose an $F$-stable Borel subgroup $B=TU$, where $U$ is its
unipotent radical and $T$ is an $F$-stable, $\theta$-stable maximal torus.
The pair $(B,T)$ exists by \cite[Lemma~3.1.1]{helminck1994symmetric}.
Write $W\coloneqq N_G(T)/T$, and let $S$ be the simple reflections
determined by $B$.
Fix an $F$-stable $W$-orbit $\fo\coloneqq W\cdot\cL$ of Kummer
local systems on $T$ (Definition~\ref{def:kummer-local-systems}).

The standard-basis parameters $\xi\in\Gamma^{H^\circ}_\fo$ consist of
an $H^\circ$-orbit in $G/B$ and an irreducible $H^\circ$-equivariant
local system on its inverse image in $G/U$, with torus monodromy
in $\fo$ (see Definitions~\ref{def:monodromic} and~\ref{def:standard-basis-parameter}).
Put $\cA\coloneqq\bZ[u^{1/2},u^{-1/2}]$.
The symmetric monodromic Hecke module
(Definition~\ref{def:relative-monodromic-hecke-module}; cf.
Henderson \cite[Chapter~3, Proposition~3.4]{henderson2003character}) is
\[
    \mathbf M^{H^\circ}_{\fo}
    \coloneqq\bigoplus_{\xi\in\Gamma^{H^\circ}_\fo}\cA m_\xi.
\]
The monodromic Hecke algebra $\mathbf H_\fo$
(Definition~\ref{def:monodromic-hecke-algebra}; see
Lusztig--Yun \cite[\S3.13]{lusztig2020endoscopy}) acts on
$\mathbf M^{H^\circ}_{\fo}$ by convolution, as constructed in
\S\ref{sec:hecke-category-and-symmetric-module}.

\begin{introtheorem}[cf. Theorem~\ref{prop:rank-one-action-formulas}]
    \label{intthm:monodromic-generator-action}
    For every $s\in S$, the action of the simple Hecke generator
    $\mathbf T_s$ on the standard basis of $\mathbf M^{H^\circ}_\fo$
    is given by the formulas in
    Theorem~\ref{prop:rank-one-action-formulas}.
    For trivial torus monodromy, these formulas coincide with the
    Lusztig--Vogan formulas \cite[Lemma~3.5]{lusztig1983singularities}.
\end{introtheorem}

Specializing $u^{1/2}$ to $-1$ gives a
$\overline{\bQ}_\ell[W]$-module $M^{H^\circ}_\fo$.
Corollary~\ref{cor:rank-filtration} gives $M^{H^\circ}_\fo$ a
$W$-stable filtration by the anti-invariant torus rank
of Definition~\ref{def:theta-degree}.
On the associated graded, each simple reflection acts by the cross
action of Definition~\ref{def:cross-action}, with sign $-1$ for
real roots and $1$ otherwise.

\subsection{Background and motivation}
\label{subsec:background-and-motivation}
Irreducible equivariant local systems on symmetric-subgroup orbits
in a flag variety index the standard basis of the Lusztig--Vogan
Hecke module. Its simple-generator action is computed in
\cite[Lemma~3.5]{lusztig1983singularities}; see also
\cite[\S7.5]{lusztig-vogan2014quasisplit}.
Henderson constructs the analogous module with prescribed torus
monodromy in \cite[Chapter~3, Proposition~3.4]{henderson2003character}.
The present paper computes the action for an entire Weyl-group orbit
of torus monodromies, including generators that change monodromy.

Our motivation comes from distinguished representations of finite
reductive groups. A representation $\pi$ of $G^F$ is
\emph{$H^F$-distinguished} if $\pi^{H^F}\ne0$.
In a subsequent paper \cite{dagan-distinction}, we use the
description of $\mathbf M^{H^\circ}_\fo$ and its specialization to
study $H^F$-distinguished irreducible representations of $G^F$.
The filtration of $M^{H^\circ}_\fo$ by anti-invariant torus rank
will yield a necessary condition for distinction in terms of
Weyl-group representations indexing Lusztig families.
\subsection{Overview of the method}
\label{subsec:overview-of-the-method}

In \S\ref{sec:standard-objects-and-simple-reflections}, we describe the standard-basis parameters by stabilizer
component-group characters and relate them through cross action
and Cayley transforms.
For each simple reflection, convolution stalks reduce to compactly
supported cohomology of punctured orbits in a projective line.
Tracking local monodromy, weights, and stabilizer characters yields
the action formulas of
\S\ref{sec:formulas-for-the-simple-hecke-operators}.

In \S\ref{sec:specialization-leading-term}, we show that
after specialization at $u^{1/2}=-1$, the cross-action term
preserves anti-invariant torus rank, while every other surviving
term strictly increases it.
These rank comparisons give the filtration and the signed
cross action on its associated graded.

\subsection{Acknowledgments}
\label{subsec:acknowledgments}
I thank Guy Shtotland for helpful conversations.
I thank my advisor, Avraham Aizenbud, for his guidance throughout this work
and for teaching me mathematics.

I was partially supported by ISF grant 1781/23 and BSF grant 2022193.

\section{Preliminaries}
\label{sec:preliminaries}
\subsection{Geometric setup}
\label{subsec:geometric-setup}
Let $G$ be a connected reductive group over $k\coloneqq\overline{\bF}_q$ with smooth
and connected center,
where $q$ is odd. Fix a Frobenius endomorphism $F:G\to G$ and an
involution $\theta$ commuting with $F$. Put $H\coloneqq G^\theta$,
and write $H^\circ$ for its identity component; both are $F$-stable.
Fix a prime $\ell\ne\opr{char}(k)$ and a compatible identification
of roots of unity of order prime to $\opr{char}(k)$ in $k$ and
$\overline{\bQ}_\ell$.

Fix an $F$-stable Borel subgroup $B$ and an $F$-stable,
$\theta$-stable maximal torus $T\subseteq B$.
Such a choice exists by
\cite[Lemma~3.1.1]{helminck1994symmetric}.
Put $U\coloneqq R_u(B)$ and $W\coloneqq N_G(T)/T$.
Let $\Phi^+$, $S$, and $\ell_W$ be the positive roots, simple
reflections, and length determined by $B$.

\subsection{Sheaves and weighted Grothendieck groups}
\label{subsec:sheaves}
For a separated scheme $X$ of finite type over $k$,
write $\cD(X)\coloneqq D_c^b(X,\overline{\bQ}_\ell)$
for the bounded constructible derived category of $\ell$-adic sheaves on $X$.
For an algebraic group $Q$ acting on $X$,
we write $\cD^Q(X)$ for the equivariant
constructible derived category as in
\cite[\S1.10.2]{lusztig2020endoscopy}.

\begin{definition}[cf. {\cite[\S\S2.1--2.2]{lusztig2020endoscopy}}]
    \label{def:kummer-local-systems}
    A \emph{Kummer local system} is a rank-one
    $\overline{\bQ}_{\ell}$-local system $\cL$ on a torus $T$
    equipped with a multiplicative isomorphism
    $m_T^*\cL\simeq\cL\boxtimes\cL$ and a trivialization of $\cL_1$,
    satisfying the associativity and unit axioms.
    Here $m_T:T\times T\to T$ is multiplication and $1\in T$ is the identity.
    We require $[n]^*\cL$ to be trivial for some $n>0$ prime to
    $\opr{char}(k)$, where $[n](t)\coloneqq t^n$.
\end{definition}
Choose $n$ as in Definition~\ref{def:kummer-local-systems} and set
$T[n]\coloneqq[n]^{-1}(1)$. There is a unique character
$\chi_{\cL}:T[n]\to\overline{\bQ}_\ell^\times$ such that
\[
    \cL\simeq([n]_*\overline{\bQ}_\ell)_{\chi_\cL}.
\]
Under the fixed identification of roots of unity, write
$\chi_{\cL}=\lambda|_{T[n]}$ with $\lambda\in X^*(T)$, determined
modulo $nX^*(T)$.
The local system $\cL$ corresponds to the rational character
class $\lambda/n+X^*(T)$ in $(X^*(T)\otimes_\bZ\bQ)/X^*(T)$.

\begin{definition}[cf. {\cite[\S2.5]{lusztig2020endoscopy}}]
    \label{def:monodromic}
    Let $a:T\times X\to X$ be a torus action and $\cL$ a Kummer
    local system. Choose $n,\chi_{\cL}$ as above, and put
    $a_n(t,x)\coloneqq a(t^n,x)$. Define
    \[
        \cD_{\cL}(X)\coloneqq
        \{K\in\cD^T(X,a_n)\mid
        z\text{ acts as }\chi_{\cL}(z)\opr{id}_K
        \text{ for all }z\in T[n]\}.
    \]
    Its objects are
    \emph{$(T,\cL)$-monodromic complexes}, with twisted equivariance
    $a^*K\simeq\cL\boxtimes K$.
\end{definition}

\begin{definition}[cf. {\cite[\S5.1]{BBD82}}]
    \label{def:weil-structure}
    Let $F_X:X\to X$ define a model $X_0$ of $X$ over $\bF_{q^r}$,
    where $r\in\bZ_{>0}$.
    A \emph{Weil structure} on $K\in\cD(X)$ is an
    isomorphism $\phi:F_X^*K\xrightarrow{\sim}K$.

    Put $\cD^{\mathrm{mix}}(X,F_X)\coloneqq
        D_m^b(X_0,\overline{\bQ}_\ell)$, the full subcategory of
    $D_c^b(X_0,\overline{\bQ}_\ell)$ whose cohomology sheaves admit
    finite filtrations with pointwise pure quotients of integral weights.
    A constructible sheaf $V$ on $X_0$ is \emph{pointwise pure of weight}
    $j\in\bZ$ if, at every closed point $x$ of degree $d$ over $\bF_{q^r}$,
    each local geometric Frobenius eigenvalue $\alpha$ on $V_{\bar x}$
    is algebraic over $\bQ$ and satisfies
    \[
        |\sigma(\alpha)|=q^{rdj/2}
        \qquad\text{for every embedding }\sigma:\bQ(\alpha)\hookrightarrow\bC.
    \]
\end{definition}

For $V\in \cD^{\mathrm{mix}}(\Spec{\bF_{q^r}})$ concentrated
in degree $0$,
write $\mathrm{Wt}_\bullet V$ for its weight filtration,
with $\opr{Gr}^{\mathrm{Wt}}_jV$ pure of weight $j$.
Put $\cA\coloneqq\bZ[u^{1/2},u^{-1/2}]$ and
$R_{\mathrm{Fr}}(q^r)\coloneqq
    K_0(\cD^{\mathrm{mix}}(\Spec{\bF_{q^r}}))$, with multiplication
induced by tensor product. Define the ring homomorphism
\[
    \opr{wt}:R_{\mathrm{Fr}}(q^r)\longrightarrow\cA,\qquad
    [V]\longmapsto\sum_{i,j}(-1)^i(-u^{1/2})^j
    \dim\opr{Gr}^{\mathrm{Wt}}_jH^i(V).
\]

\begin{definition}[cf. {\cite[\S\S1.10.1--1.10.2]{lusztig2020endoscopy}}]
    \label{def:tate-twists-and-normalization}

    Let $\mathcal C^{\mathrm{mix}}$ be a triangulated category with an
    action of $\cD^{\mathrm{mix}}(\Spec{\bF_{q^r}})$ exact in each variable.
    Its \emph{weighted Grothendieck group} and the \emph{weighted class}
    of $V\in\mathcal C^{\mathrm{mix}}$ are
    \[
        K_0^{\mathrm{wt}}(\mathcal C^{\mathrm{mix}})
        \coloneqq\cA\otimes_{R_{\mathrm{Fr}}(q^r),\opr{wt}}
        K_0(\mathcal C^{\mathrm{mix}}),\qquad
        \chi_u(V)\coloneqq1\otimes[V].
    \]

\end{definition}
For the equivariant monodromic categories
$\mathcal C^{\mathrm{mix}}=\cD^{\mathrm{mix},Q}_\cL(X,F_X)$ used below,
the commuting $Q$- and $T$-actions have finitely many $Q\times T$-orbits.
Assume that $F_X$ preserves every orbit and $F_X^*$ fixes every isomorphism
class of irreducible $Q$-equivariant $(T,\cL)$-monodromic local systems
on those orbits.
The weight-polynomial $\chi_u$ description of
\cite[pp.~27--28]{henderson2003character} identifies the weighted group
with the free $\cA$-module on simple perverse objects,
equipped with pure weight-zero Weil structures.
Base change preserves these classes and their weights. Hence base change
to $\bF_{q^{rm}}$, with $m\geq1$, induces an $\cA$-linear
isomorphism
\[
    K_0^{\mathrm{wt}}\!\left(\cD^{\mathrm{mix},Q}_\cL(X,F_X)\right)
    \xrightarrow{\sim}
    K_0^{\mathrm{wt}}\!\left(\cD^{\mathrm{mix},Q}_\cL(X,F_X^m)\right).
\]
For an irreducible equivariant monodromic local system $\mathscr E$
on one orbit, two compatible pure weight-zero Weil structures
$\phi,\phi'$ differ by a scalar of weight zero. Thus
$\chi_u(\mathscr E,\phi)=\chi_u(\mathscr E,\phi')$ in the weighted
Grothendieck group of that orbit.

For a finite group $L$ acting trivially on $\Spec{\bF_{q^r}}$,
the natural identification
\[
    K_0^{\mathrm{wt}}\bigl(\cD^{\mathrm{mix},L}(\Spec{\bF_{q^r}})\bigr)
    \simeq\cA\otimes_\bZ K_0(\opr{Rep}_{\overline{\bQ}_\ell}(L))
\]
expresses the weighted class of a mixed complex $V$ with an $L$-action
commuting with Frobenius as the \emph{representation-valued weight polynomial}
\[
    \chi_u(V)=\sum_{i,j}(-1)^i(-u^{1/2})^j
    [\opr{Gr}^{\mathrm{Wt}}_jH^i(V)].
\]
For $L=1$, this is the scalar weight polynomial $\opr{wt}([V])$.

Fix a square root $q^{1/2}\in\overline{\bQ}_\ell$ and let $\overline{\bQ}_\ell(1/2)$ be
the one-dimensional sheaf on $\Spec{\bF_{q^r}}$ with Frobenius
eigenvalue $(q^{1/2})^{-r}$. For $n\in\bZ$, put
\[
    V(n/2)\coloneqq V\otimes\overline{\bQ}_\ell(1/2)^{\otimes n},
    \qquad V\langle n\rangle\coloneqq V[n](n/2).
\]
Here $[n]$ is the cohomological shift and $V(1)$ is the usual Tate twist.

An exact functor $F:\mathcal C^{\mathrm{mix}}\to(\mathcal C')^{\mathrm{mix}}$
compatible with the scalar action induces an $\cA$-linear map $[F]$
on weighted Grothendieck groups, with $[F](\chi_u(V))=\chi_u(FV)$.

Triangle additivity and the weight of the half Tate twist give
\begin{equation}
    \label{eq:weighted-shift-conventions}
    \begin{split}
        \chi_u(V[1])              & =-\chi_u(V),\qquad
        \chi_u(V(1/2))=-u^{-1/2}\chi_u(V),              \\
        \chi_u(V\langle n\rangle) & =u^{-n/2}\chi_u(V).
    \end{split}
\end{equation}
Evaluation at $u^{1/2}=-1$ sends a Frobenius vector space to its
dimension, and a complex on a point to its Euler characteristic.

\subsection{\'Etale cohomology}
\label{subsec:etale-cohomology}
Let $K$ be an algebraic group over $k$, possibly disconnected, acting
on a $k$-variety $D$.
Let $\mathscr F$ be a $K$-equivariant rank-one local system on $D$.
Assume that $K$, $D$, the action, and the equivariant local system
are defined over a finite subfield of $k$.
Equip $\mathscr F$ with a Weil structure pure of weight zero.

The equivariant structure induces a $K$-action on each
$H_c^i(D,\mathscr F)$. All cohomology is taken over $k$, and the
decompositions below are isomorphisms of $K$-representations.
We record Frobenius weights separately.
A Tate twist changes the weights but leaves the $K$-action unchanged.

\begin{lemma}
    \label{lem:finite-domain-cohomology}
    If $D$ is finite, then
    \[
        H_c^0(D,\mathscr F)\cong\bigoplus_{x\in D}\mathscr F_x.
    \]
    The group $H_c^0(D,\mathscr F)$ has weight zero, and
    $H_c^i(D,\mathscr F)=0$ for $i\ne0$.
\end{lemma}
\begin{proof}
    The decomposition follows by writing $D$ as a disjoint union of
    geometric points. An element $g\in K$ sends the summand
    $\mathscr F_x$ to $\mathscr F_{gx}$ through the equivariant structure.
    Purity of $\mathscr F$ gives the weight assertion.
    If $D=\emptyset$, all cohomology groups vanish.
\end{proof}

For the rest of this subsection, assume that $K$ acts on $\bP_1$.
Let $Z\subset\bP_1$ be a nonempty finite set preserved by $K$, and
choose $z\in Z$ fixed by  $K$.
Thus $K$ may permute the points of $Z\setminus\{z\}$.
Put $D\coloneqq\bP_1\setminus Z$, with the restricted $K$-action.

Assume also that $\mathscr F^{\otimes N}$ is geometrically constant
for some positive integer $N$ prime to $\opr{char}(k)$.
Equivalently, the geometric monodromy character of $\mathscr F$
has finite image of order dividing $N$.
In particular, $\mathscr F$ is tame at every point of $Z$.
The local monodromy of $\mathscr F$ at $c\in Z$ is trivial exactly when
$\mathscr F$ extends as a local system across $c$.

Since $D$ is a smooth connected nonproper curve, its compactly
supported cohomology with coefficients in $\mathscr F$ can be
nonzero only in degrees one and two.
For geometrically constant $\mathscr F$, write $\rho$ for the
$K$-character on the one-dimensional space $H^0(D,\mathscr F)$.
Localization and Poincar\'e duality
\cite[Proposition~18.3(a) and Theorem~24.1]{milne2013etale}
give weights zero and two for $H_c^1(D,\mathscr F)$ and
$H_c^2(D,\mathscr F)$, respectively.

\begin{lemma}
    \label{lem:global-section-character-from-stabilizer}
    Suppose in addition that $\mathscr F$ is geometrically constant
    and that the identity component $K^\circ$ acts transitively on $D$.
    For $x\in D$, the natural map
    \[
        \opr{Stab}_K(x)\longrightarrow\pi_0(K)
    \]
    is surjective. Evaluation at $x$ identifies the restriction of
    $\rho$ to $\opr{Stab}_K(x)$ with the character on $\mathscr F_x$.
    The character on $\mathscr F_x$ determines $\rho$ uniquely.
\end{lemma}
\begin{proof}
    Transitivity gives $K=K^\circ\opr{Stab}_K(x)$,
    proving surjectivity.
    Evaluation at $x$ identifies $H^0(D,\mathscr F)$ with
    $\mathscr F_x$ equivariantly for $\opr{Stab}_K(x)$.
    The action on global sections is locally constant on $K$,
    so $\rho$ factors through $\pi_0(K)$.
    Its restriction to $\opr{Stab}_K(x)$ determines $\rho$.
\end{proof}

\begin{lemma}
    \label{lem:affine-line-cohomology}
    If $|Z|=1$, then $D\cong\bA_1$, the local system $\mathscr F$
    is geometrically constant, and
    \[
        H_c^1(D,\mathscr F)=0,\qquad
        H_c^2(D,\mathscr F)\cong\rho(-1).
    \]
\end{lemma}
\begin{proof}
    The Kummer calculation in \cite[Example~16.3]{milne2013etale}
    shows that $\mathscr F$ is geometrically constant, since its
    geometric monodromy has order prime to $\opr{char}(k)$.
    Apply Poincar\'e duality \cite[Theorem~24.1]{milne2013etale}
    to the affine-line calculation in \cite[Example~16.3]{milne2013etale}.
    The trace map is $K$-invariant.
\end{proof}

\begin{lemma}
    \label{lem:multiplicative-group-cohomology}
    Suppose $|Z|=2$, so $D\cong\bG_m$.
    If $\mathscr F$ is geometrically constant, then
    \[
        H_c^1(D,\mathscr F)\cong\rho,\qquad
        H_c^2(D,\mathscr F)\cong\rho(-1).
    \]
    Otherwise, all compactly supported cohomology groups vanish.
\end{lemma}
\begin{proof}
    Since $K$ preserves $Z$ and fixes $z$, it fixes both points of $Z$.
    If $\mathscr F$ is geometrically constant, extend it across the
    point of $Z\setminus\{z\}$ to $\bP_1\setminus\{z\}\cong\bA_1$.
    The extension of $\mathscr F$ is $K$-equivariant. Apply localization
    \cite[Proposition~18.3(a)]{milne2013etale} with
    Lemmas~\ref{lem:finite-domain-cohomology} and
    \ref{lem:affine-line-cohomology}.

    If $\mathscr F$ is geometrically nonconstant, duality and the
    tame Euler-characteristic calculation in
    \cite[proof of Lemma~2.2]{katz2013squarefree}
    give $H_c^2(D,\mathscr F)=0$ and
    $\dim H_c^1(D,\mathscr F)=-\chi_c(D)=0$.
\end{proof}

When $|Z|=3$, choose a coordinate $t$ on $\bP_1$ sending $z$ to
$\infty$ and the other two points of $Z$ to $0$ and $1$.
Every element of $K$ fixes $\infty$ and preserves $\{0,1\}$.
An automorphism of $\bP_1$ fixing three distinct points is the identity.
Hence each element of $K$ acts by either $t\mapsto t$ or $t\mapsto1-t$.
In particular, the image of $K$ in $\opr{Aut}(\bP_1)$ has order at most two.

Let $\delta$ be the sign character of the $K$-action on $Z\setminus\{z\}$.
Thus $\delta(g)=1$ when $g$ acts trivially on $D$, and $\delta(g)=-1$
when $g$ acts by $t\mapsto1-t$.
Since $\opr{char}(k)\ne2$, the point $x$ with $t(x)=1/2$ belongs to $D$
and is fixed by every element of $K$.
Write $\rho_x$ for the character on $\mathscr F_x$.
For geometrically constant $\mathscr F$, evaluation at $x$ identifies
$\rho_x$ with $\rho$.

\begin{lemma}
    \label{lem:action-on-cohomology}
    If $|Z|=3$, the compactly supported cohomology is
    \[
        \begin{array}{c|c|c}
            \mathscr F                       & H_c^1(D,\mathscr F)  & H_c^2(D,\mathscr F) \\ \hline
            \text{geometrically constant}    & \rho\oplus\rho\delta & \rho(-1)            \\
            \text{geometrically nonconstant} & \rho_x\delta         & 0
        \end{array}
    \]
    If $\mathscr F$ is geometrically nonconstant, then $H_c^1(D,\mathscr F)$
    has weight zero when the monodromy at some $c\in Z$ is trivial,
    and weight one otherwise.
\end{lemma}
\begin{proof}
    If $\mathscr F$ is geometrically constant, localization in
    $\bP_1\setminus\{z\}$ and Lemmas~\ref{lem:finite-domain-cohomology}
    and \ref{lem:affine-line-cohomology} give the first row of the table:
    the permutation representation on $Z\setminus\{z\}$ is $1\oplus\delta$.

    Suppose now that $\mathscr F$ is geometrically nonconstant.
    In the coordinate $t$, it is geometrically a tensor product of
    Kummer sheaves associated with $t$ and $t-1$.
    Apply \cite[Lemma~2.2(1)--(2) and its proof]{katz2013squarefree}
    with $n=2$ and $\chi_1,\chi_2$ the monodromies at $0,1$.
    After a finite field extension, $\mathscr F$ differs from Katz's sheaf
    by a geometrically constant twist of weight zero.
    Thus $\dim H_c^1(D,\mathscr F)=1$, $H_c^2(D,\mathscr F)=0$,
    and the weight assertion follows.

    To determine the $K$-character of $H_c^1(D,\mathscr F)$, fix $g\in K$.
    If $\delta(g)=1$, then $g$ acts trivially on $D$.
    Since $D$ is connected, $g$ acts on $\mathscr F$ by the scalar
    $\rho_x(g)$, and hence acts by the same scalar on cohomology.

    Suppose $\delta(g)=-1$, and put $\sigma(t)\coloneqq1-t$.
    The isomorphism $\sigma^*\mathscr F\cong\mathscr F$ forces the
    monodromies at $0$ and $1$ to agree.
    The Kummer description gives $\mathscr F\cong\pi^*\mathscr L$
    for a geometrically nonconstant Kummer local system $\mathscr L$
    on $\bG_m$ of order dividing $N$, where
    \[
        \pi:D\longrightarrow\bG_m,\qquad \pi(t)\coloneqq t(t-1),
    \]
    is the finite quotient by $\sigma$.
    The canonical $\sigma$-action on $\pi^*\mathscr L$ is the identity
    on $\mathscr F_x$.
    Hence the coefficient map of $g$ is $\rho_x(g)$ times the
    canonical $\sigma$-action.

    Since $(\pi_*\overline{\bQ}_\ell)^\sigma\cong\overline{\bQ}_\ell$,
    the projection formula \cite[Proposition~22.5]{milne2013etale} gives
    \[
        H_c^1(D,\mathscr F)^\sigma
        \cong H_c^1(\bG_m,\mathscr L)=0,
    \]
    where the vanishing follows from
    Lemma~\ref{lem:multiplicative-group-cohomology} with trivial acting group.
    Thus $\sigma$ acts by $-1$ on the one-dimensional $H_c^1(D,\mathscr F)$,
    and $g$ acts by $-\rho_x(g)=\rho_x(g)\delta(g)$.
\end{proof}

\subsection{Relative Weyl groups and the Kummer exponent}
\label{subsec:relative-weyl-groups}
\begin{definition}[cf. {\cite[\S2.4]{lusztig2020endoscopy}}]
    \label{def:relative-weyl-group}
    For a Kummer local system $\cL$ on $T$, define
    \[
        W'_\cL\coloneqq\{w\in W\mid w\cL\cong\cL\},
        \qquad
        \Phi_\cL\coloneqq\{\alpha\in\Phi(G,T)\mid
        (\alpha^\vee)^*\cL\cong\overline{\bQ}_\ell\},
    \]
    where $w\cL\coloneqq(w^{-1})^*\cL$.
    The set $\Phi_\cL$ is a root subsystem.
    Set
    \[
        \Phi^+_\cL\coloneqq\Phi_\cL\cap\Phi^+,\qquad
        W^\circ_\cL\coloneqq\langle s_\alpha\mid\alpha\in\Phi_\cL\rangle.
    \]
    Write $S_\cL$ and $\ell_\cL$ for the simple reflections and
    length of $W^\circ_\cL$ determined by $\Phi^+_\cL$.
\end{definition}
In general, $S_\cL$ is not a subset of $S$, and $\ell_\cL$ is not
the restriction of $\ell_W$.

\begin{lemma}[cf. {\cite[\S9.1]{lusztig2020endoscopy}} and
            {\cite[Theorem~5.13]{deligneLusztig1976}}]
    \label{lem:connected-center-character-weyl-group}
    If $G$ has connected center, then
    \[
        W'_\cL=W^\circ_\cL=W(\Phi_\cL),
    \]
    for every Kummer local system $\cL$ on $T$.
\end{lemma}
Under our standing connected-center assumption, write
$W_\cL\coloneqq W'_\cL=W^\circ_\cL$ from now on.

\begin{remark}[cf. {\cite[\S2.2]{lusztig2020endoscopy}}]
    \label{rem:uniform-kummer-exponent}
    Throughout, put $p\coloneqq\opr{char}(k)$ and
    $N\coloneqq |G^F|_{p'}$, the largest divisor of $|G^F|$
    prime to $p$.
    If $\cL$ is a Kummer local system on $T$ with $F$-stable
    $W$-orbit, then
    \[
        [N]^*\cL\cong\overline{\bQ}_\ell.
    \]
\end{remark}

\section{The Hecke category and the symmetric module}
\label{sec:hecke-category-and-symmetric-module}
Let $B=TU$ be as in \S\ref{subsec:geometric-setup}, and fix the
$F$-stable orbit $\fo\coloneqq W\cdot\cL$. Set
\[
    \cB\coloneqq G/B,\qquad \widetilde{\cB}\coloneqq G/U,\qquad
    \pi:\widetilde{\cB}\longrightarrow\cB.
\]

\subsubsection{Monodromic categories and convolution}
\label{subsubsec:monodromic-categories-and-convolution}
As in \cite[pp.~27--28]{henderson2003character}, choose an auxiliary
field $\bF_{q^r}$ over which $T$ splits and all elements of $\fo$
are defined. Enlarge $r$ so that every orbit and every irreducible
equivariant local system used below is fixed by $F^{r*}$.
Choose pure weight-zero structures on these local systems.
Such a common choice exists because the orbit and local-system
parameter sets are finite. We use the mixed categories of
Definition~\ref{def:weil-structure} over this auxiliary field.

\begin{definition}[cf. {\cite[\S2.7]{lusztig2020endoscopy}}]
    \label{def:hecke-category}
    For $\cL_1,\cL_2\in\fo$, the \emph{mixed Hecke category with monodromies $(\cL_1,\cL_2)$} is
    \[
        \cH_{\cL_1,\cL_2}
        \coloneqq\cD^{G,\opr{mix}}_{\cL_1\boxtimes\cL_2^{-1}}(\widetilde{\cB}^2),
    \]
    where the action is
    \[
        (g,t_1,t_2)\cdot(g_1U,g_2U)
        =(gg_1t_1^{-1}U,gg_2t_2^{-1}U).
    \]
    Write $\cH_{\cL'}\coloneqq\cH_{\cL',\cL'}$. The \emph{mixed Hecke category for the orbit $\fo$} is
    \[
        \cH_{\fo}
        \coloneqq\bigoplus_{\cL_1,\cL_2\in\fo}\cH_{\cL_1,\cL_2}.
    \]
\end{definition}

\begin{definition}
    \label{def:relative-module-category}
    For a Kummer local system $\cL'$ on $T$, the
    \emph{mixed symmetric Hecke module category with monodromy $\cL'$} and the
    \emph{mixed symmetric Hecke module category for the orbit $\fo$} are
    \[
        \cM^{H^\circ}_{\cL'}\coloneqq\cD^{H^\circ,\opr{mix}}_{\cL'}(\widetilde{\cB}),\qquad
        \cM^{H^\circ}_{\fo}
        \coloneqq\bigoplus_{\cL'\in\fo}\cM^{H^\circ}_{\cL'},
    \]
    respectively. The action defining $\cM^{H^\circ}_{\cL'}$ is
    $(h,t)\cdot gU=hgt^{-1}U$, with $(h,t)\in H^\circ\times T$.
    The adjective ``symmetric'' refers to the symmetric space $G/H^\circ$.
\end{definition}

Let $p_{ij}:\widetilde{\cB}^3\to\widetilde{\cB}^2$ be the projections.
Factor $p_{13}$ as
\[
    \widetilde{\cB}^3
    \xrightarrow{\ q=\opr{id}\times\pi\times\opr{id}\ }
    \widetilde{\cB}\times\cB\times\widetilde{\cB}
    \xrightarrow{\ \overline p_{13}\ }\widetilde{\cB}^2.
\]
For $K\in\cH_{\cL_1,\cL_2}$ and $K'\in\cH_{\cL_2,\cL_3}$,
the monodromies of $p_{12}^*K\otimes p_{23}^*K'$ cancel on the middle copy of $\widetilde{\cB}$.
The resulting middle $T$-equivariance allows descent along $q$.
Write $q_\flat$ for the inverse of the pullback $q^*$. The projection $\overline p_{13}$ is proper,
with fiber $\cB$. Recall the normalized shift
$K\langle n\rangle\coloneqq K[n](n/2)$ of
Definition~\ref{def:tate-twists-and-normalization}, and put
\[
    d_0\coloneqq\dim\widetilde{\cB}+\dim T.
\]

\begin{definition}
    \label{def:hecke-convolution}
    For $K\in\cH_{\cL_1,\cL_2}$ and $K'\in\cH_{\cL_3,\cL_4}$, define their \emph{convolution} by
    \[
        K\star K'
        \coloneqq
        \begin{cases}
            \overline p_{13!}\,q_\flat
            (p_{12}^*K\otimes p_{23}^*K')\langle-d_0\rangle
            \in\cH_{\cL_1,\cL_4} & \text{if } \cL_2 \cong \cL_3 \\
            0                    & \text{otherwise.}
        \end{cases}
    \]
\end{definition}

For the action of $\cH_\fo$ on $\cM^{H^\circ}_\fo$, let $p_2:\widetilde{\cB}^2\to\widetilde{\cB}$
be the second projection and consider
\[
    \widetilde{\cB}^2
    \xrightarrow{\ q_2=\opr{id}\times\pi\ }
    \widetilde{\cB}\times\cB
    \xrightarrow{\ \overline p_1\ }\widetilde{\cB}.
\]
If $A\in\cM^{H^\circ}_{\cL_2}$ and $K\in\cH_{\cL_1,\cL_2}$, then $K\otimes p_2^*A$ descends along $q_2$. Write $q_{2\flat}$ for the inverse of the pullback $q_2^*$.
The projection $\overline p_1$ is proper.

\begin{definition}
    \label{def:relative-convolution}
    For $K\in\cH_{\cL_1,\cL_2}$ and $A\in\cM^{H^\circ}_{\cL_3}$, define their \emph{convolution} by
    \[
        K\star A
        \coloneqq
        \begin{cases}
            \overline p_{1!}\,q_{2\flat}
            (K\otimes p_2^*A)\langle-d_0\rangle
            \in\cM^{H^\circ}_{\cL_1} & \text{if } \cL_2 \cong \cL_3 \\
            0                        & \text{otherwise.}
        \end{cases}
    \]
\end{definition}

The following lemma is standard.
\begin{lemma}[cf. {\cite[\S3.1]{lusztig2020endoscopy}} and
            {\cite[pp.~26--28]{henderson2003character}}]
    \label{lem:convolution-associativity}
    Convolution commutes with forgetting the mixed structure.
    It preserves distinguished triangles and shifts in each variable.
    There are compatible associativity isomorphisms
    \[
        \begin{gathered}
            (K_{12}\star K_{23})\star K_{34}        \simeq K_{12}\star(K_{23}\star K_{34}), \\
            (K_{12}\star K_{23})\star A_3          \simeq K_{12}\star(K_{23}\star A_3)    \\
            (K_{12}\otimes V)\star A_2  \simeq (K_{12}\star A_2)\otimes V
            \simeq K_{12}\star(A_2\otimes V),                                \\
            (K_{12}\otimes V)\star K_{23} \simeq (K_{12}\star K_{23})\otimes V
            \simeq K_{12}\star(K_{23}\otimes V),\\
            F^*(K\star K')\simeq F^*K\star F^*K',\qquad
            F^*(K\star A)\simeq F^*K\star F^*A,
        \end{gathered}
    \]
    where $K_{i,i+1} \in \cH_{\cL_i,\cL_{i+1}}$ and $A_i\in\cM^{H^\circ}_{\cL_i}$ for $i \leq 3$, and $V\in \cD^{\opr{mix}}(\opr{pt})$.
\end{lemma}

\subsubsection{Standard objects}
\label{subsubsec:standard-objects}
For each $w\in W$, choose a representative $\dot w$, with $\dot1=1$.
The associated $G\times T\times T$-orbit in $\widetilde{\cB}^2$ is
\[
    \widetilde{\cO}_w
    \coloneqq\{(g_1U,g_2U)\mid g_1^{-1}g_2\in B\dot wB\}\subseteq \widetilde{\cB}^2,
\]
and $j_w:\widetilde{\cO}_w\hookrightarrow\widetilde{\cB}^2$ denotes its inclusion.
The orbit $\widetilde{\cO}_w$ has dimension $d_0+\ell_W(w)$.
By \cite[Prop.~2.2]{henderson2003character}, for any $\cL'\in\fo$ there is a unique irreducible
equivariant local system $\mathscr E_{w,\cL'}$ on $\widetilde{\cO}_w$ with monodromy $w\cL'\boxtimes(\cL')^{-1}$. Give $\mathscr E_{w,\cL'}$ a pure $\bF_{q^r}$-structure of weight $0$.

\begin{definition}[cf. {\cite[\S2.11]{lusztig2020endoscopy}}]
    \label{def:hecke-standard-objects}
    The \emph{standard object} in $\cH_\fo$ indexed by
    $(w,\cL')$ is
    \[
        \Delta_{w,\cL'}\coloneqq j_{w!}
        \bigl(\mathscr E_{w,\cL'}\langle d_0+\ell_W(w)\rangle\bigr)\in \cH_{w\cL', \cL'}\subseteq \cH_\fo.
    \]
\end{definition}

\begin{definition}
    \label{def:standard-basis-parameter}
    For a Kummer local system $\cL'$ on $T$, let $\Gamma^{H^\circ}_{\cL'}$ be the set of pairs $\xi=(\cO_\xi,\mathscr E_\xi)$,
    where
    \begin{itemize}
        \item $\cO_\xi$ is an $H^\circ$-orbit in $\cB$;
        \item $\mathscr E_\xi$ is a representative of an isomorphism class of irreducible $H^\circ$-equivariant $(T,\cL')$-monodromic local systems on $\widetilde{\cO}_\xi\coloneqq\pi^{-1}(\cO_\xi)$.
    \end{itemize}
    Put
    \[
        \Gamma^{H^\circ}_{\fo}\coloneqq\coprod_{\cL'\in\fo}
        \Gamma^{H^\circ}_{\cL'}.
    \]
    For $\xi\in\Gamma^{H^\circ}_{\cL'}$, write
    $\cL_\xi\coloneqq\cL'$ for its torus monodromy.
\end{definition}
We give each $\mathscr E_\xi$ a pure $\bF_{q^r}$-structure of weight $0$.
For $\xi\in\Gamma^{H^\circ}_{\fo}$, write
$j_\xi:\widetilde{\cO}_\xi\hookrightarrow\widetilde{\cB}$ and
$d_\xi\coloneqq\dim\widetilde{\cO}_\xi=\dim\cO_\xi+\dim T$.
\begin{definition}
    \label{def:relative-standard-objects}
    The \emph{standard object} in $\cM^{H^\circ}_\fo$ is
    \[
        \Delta_\xi\coloneqq j_{\xi!}(\mathscr E_\xi\langle d_\xi\rangle).
    \]
\end{definition}

\subsubsection{Standard bases and the generic Hecke algebra}
\label{subsubsec:the-generic-module-and-frobenius}
We use the weight-polynomial description of
\cite[pp.~27--28]{henderson2003character} for the weighted Grothendieck
groups of Definition~\ref{def:tate-twists-and-normalization}.
The standard classes and intersection-cohomology classes differ by a
triangular change of basis with invertible diagonal.
The next definition follows \cite[\S3.13]{lusztig2020endoscopy}.
\begin{definition}[{\cite[1.4]{lusztig2017conjugacy}}]
    \label{def:monodromic-hecke-algebra}
    The \emph{monodromic Hecke algebra}, denoted by $\mathbf H_{\fo}$, is the unital associative $\cA$-algebra with generators $\mathbf T_w$ for $w\in W$ and $1_{\cL'}$ $(\cL' \in \fo)$ and relations
    \[
        \renewcommand{\arraystretch}{1.6}
        \begin{array}{@{}ll@{}}
            1_{\mathcal L}1_{\mathcal L'}
            = \delta_{\mathcal L,\mathcal L'}1_{\mathcal L},
             & \text{for } \mathcal L,\mathcal L' \in \fo; \\

            \mathbf T_w\mathbf T_{w'} = \mathbf T_{ww'},
             & \text{if } w,w' \in W
            \text{ and } \ell_W(ww')=\ell_W(w)+\ell_W(w'); \\

            \mathbf T_w1_{\mathcal L} = 1_{w\mathcal L}\mathbf T_w,
             & \text{for } w\in W,\ \mathcal L\in\fo;      \\

            \displaystyle
            \mathbf T_s^2 = u\mathbf T_1+(u-1)
            \sum_{\mathcal L;\,s\in W_{\mathcal L}}
            \mathbf T_s1_{\mathcal L},
             & \text{for simple reflections } s\in W;      \\[6pt]

            \displaystyle
            \mathbf T_1 = 1=\sum_{\mathcal L\in\fo}1_{\mathcal L}.
             &
        \end{array}
    \]
    Set $\mathbf{H}_{\cL_1,\cL_2} \coloneqq 1_{\cL_1} \mathbf{H}_{\fo} 1_{\cL_2}$ ($\cL_1,\cL_2 \in \fo$), which is generated by $\{1_{\cL_1} \mathbf T_w 1_{\cL_2}\,|\, w\in W,\, w\cL_2 \cong \cL_1   \}$. Write $\mathbf{H}_{\cL'}\coloneqq \mathbf{H}_{\cL',\cL'}$.
\end{definition}
By \cite[\S3.13]{lusztig2020endoscopy},
$\{\mathbf T_w1_{\cL'};(w,\cL')\in W\times\fo\}$ is an $\cA$-basis of
$\mathbf H_\fo$, and the elements with $w\in W_{\cL'}$ form a
basis of $\mathbf H_{\cL'}$. The map $\chi_u$ induces isomorphisms
\[
    K_0^{\mathrm{wt}}(\cH_{\fo})\cong\mathbf{H}_\fo, \qquad
    K_0^{\mathrm{wt}}(\cH_{\cL_1,\cL_2})\cong\mathbf{H}_{\cL_1,\cL_2},
\]
such that $\chi_u(\Delta_{w,\cL'})=u^{-\ell_W(w)/2}\mathbf T_w1_{\cL'}$.
The only nonzero stalk cohomology of $\Delta_{w,\cL'}$ at
$(U,\dot wU)$ has weight $-d_0-\ell_W(w)$.
This explains the normalization
$\langle-d_0\rangle$ in Definitions~\ref{def:hecke-convolution} and \ref{def:relative-convolution}.
\begin{definition}[cf. {\cite[\S3]{henderson2003character}}]
    \label{def:relative-monodromic-hecke-module}
    For a Kummer local system $\cL'$ on $T$, let
    $\mathbf M^{H^\circ}_{\cL'}$ be the free $\cA$-module with
    basis $m_\xi$ for $\xi\in\Gamma^{H^\circ}_{\cL'}$.
    The \emph{symmetric monodromic Hecke module} for $\fo$ is
    \[
        \mathbf M^{H^\circ}_{\fo}\coloneqq
        \bigoplus_{\cL'\in\fo}\mathbf M^{H^\circ}_{\cL'}.
    \]
\end{definition}
Similarly, \cite[p.~28]{henderson2003character} gives isomorphisms
\[
    K_0^{\mathrm{wt}}(\cM^{H^\circ}_{\fo})\cong\mathbf{M}^{H^\circ}_\fo, \qquad
    K_0^{\mathrm{wt}}(\cM^{H^\circ}_{\cL'})\cong\mathbf{M}^{H^\circ}_{\cL'},
\]
such that $\chi_u(\Delta_{\xi})=u^{-d_\xi/2}m_\xi$.
By \cite[\S3.13]{lusztig2020endoscopy} and
\cite[Theorem~3.2 and Proposition~3.4]{henderson2003character},
convolution gives a left $\mathbf H_{\fo}$-module
$\mathbf M^{H^\circ}_{\fo}$ with
\[
    \chi_u(K)\chi_u(A)=\chi_u(K\star A), \qquad
    A\in \cM^{H^\circ}_{\fo},\,
    K\in \cH_{\fo}.
\]
For each $\cL'\in\fo$, one has
$1_{\cL'}\mathbf M^{H^\circ}_{\fo}=\mathbf M^{H^\circ}_{\cL'}$.
Lemma~\ref{lem:convolution-associativity} also gives a compatible
action of Frobenius.
Therefore we obtain
\begin{equation}
    \label{eq:convolution-formula-for-hecke-action}
    \mathbf T_s m_\xi
    =
    \chi_u\!\left(
    K_{s,\xi}
    \right),\qquad
    K_{s,\xi}\coloneqq
    \overline p_{1!}\,q_{2\flat}
    \bigl((j_{s!}\mathscr E_{s,\cL_\xi})
    \otimes p_2^*(j_{\xi!}\mathscr E_\xi)\bigr),
    \qquad
    s\in S,\quad \xi\in\Gamma^{H^\circ}_{\fo}.
\end{equation}

Set $\epsilon:\cA\longrightarrow\overline{\bQ}_\ell$ by $\epsilon(u^{1/2})=-1$. We also use $\epsilon$ to denote the map $\epsilon:N \to N\otimes_{\cA,\epsilon}\overline{\bQ}_\ell$ for any $\cA$-module $N$ defined by $x\mapsto x\otimes 1$. Then $\epsilon$ induces an isomorphism
\[
    \overline{\bQ}_\ell[W\ltimes \fo]
    \simeq
    \mathbf H_{\fo}
    \otimes_{\cA,\epsilon}
    \overline{\bQ}_\ell,\qquad
    \overline{\bQ}_\ell[W_{\cL'}]
    \simeq
    \mathbf H_{\cL'}
    \otimes_{\cA,\epsilon}
    \overline{\bQ}_\ell,
\]
where $\overline{\bQ}_\ell[W\ltimes \fo]$ is the action groupoid algebra. Convolution then gives an action of $W$ and $W_{\cL'}$ on
\[
    M^{H^\circ}_{\fo} \coloneqq\overline{\bQ}_\ell\otimes_{\cA,\epsilon}\mathbf M^{H^\circ}_{\fo},\qquad
    M^{H^\circ}_{\cL'} \coloneqq\overline{\bQ}_\ell\otimes_{\cA,\epsilon}\mathbf M^{H^\circ}_{\cL'},
\]
respectively.
We also write $m_\xi$ for the image of the generic standard basis
vector under specialization.

\section{Standard objects and simple reflections}
\label{sec:standard-objects-and-simple-reflections}
\subsection{Classification of standard objects}
\label{subsec:classification-of-standard-objects}
We describe the standard-basis parameter set $\Gamma^{H^\circ}_{\fo}$
of Definition~\ref{def:standard-basis-parameter} in terms of component-group characters.
By \cite[Proposition~6.8]{helminck1993rationality},
for any $H^\circ$-orbit $\cO$ in $\cB$ there is an element $f\in G$,
unique up to left $H^\circ$-multiplication and right $T$-multiplication,
with $f^{-1}\theta(f)\in N_G(T)$ and $\cO=H^\circ fB/B$.
Let $\mathscr V$ contain one such representative $f$ for each orbit.

Assume that $\fo$ is $F$-stable.
By Remark~\ref{rem:uniform-kummer-exponent}, $[N]^*\cL'$ is
trivial for every $\cL'\in\fo$.
For any $f\in \mathscr V$, set
\begin{gather*}
    \theta_{f} \coloneqq \opr{Int}_{f^{-1}}\circ \theta \circ \opr{Int}_{f},\qquad
    H_f \coloneqq f^{-1}Hf=G^{\theta_f},\qquad
    T_{f,N} \coloneqq [N]^{-1}(T\cap H_f^\circ).
\end{gather*}
Note that $T$ is $\theta_f$-stable if and only if $f^{-1}\theta(f)\in N_G(T)$.
For a finite abelian group $A$, write
$\widehat A\coloneqq\Hom(A,\overline{\bQ}_\ell^\times)$.
We also use $\mu\in\widehat{\pi_0(T_{f,N})}$ to denote its
pullback to $T_{f,N}$ along the component map.
For $f\in\mathscr V$ and
$\mu\in\widehat{\pi_0(T_{f,N})}$, put
\[
    \cL_{f,\mu}\coloneqq
    \bigl([N]_*\overline{\bQ}_\ell\bigr)_{\mu|_{T[N]}}.
\]
\begin{lemma}
    \label{lem:combi-parameter-rel-standard-obj}
    We have a bijection
    \[
        \Gamma^{H^\circ}_{\fo}
        \longleftrightarrow
        \Big\{
        (f,\mu) \,|\,
        f\in \mathscr{V},\,
        \mu \in \widehat{\pi_0(T_{f,N})},\,
        \cL_{f,\mu} \in \fo
        \Big\}.
    \]
    Under this bijection, $(f,\mu)\in\Gamma^{H^\circ}_{\cL'}$
    if and only if $\cL_{f,\mu}\cong\cL'$.
    In particular, $\cL_\xi\cong\cL_{f,\mu}$ for $\xi=(f,\mu)$.
\end{lemma}
\begin{proof}
    Fix $f\in\mathscr V$ and put
    $\cO\coloneqq H^\circ fB/B$.
    Left multiplication by $f^{-1}$ identifies
    $\pi^{-1}(\cO)$ with $H_f^\circ B/U$ and conjugates $H^\circ$ to $H_f^\circ$.
    We classify irreducible $H_f^\circ\times T$-equivariant local systems
    on $H_f^\circ B/U$ for the action
    \[
        (a,t)\cdot xU\coloneqq ax t^{-N}U.
    \]

    Since $[N]:T\to T$ is surjective, this action is transitive.
    Its stabilizer at $U$ is
    \[
        S_f\coloneqq\opr{Stab}_{H_f^\circ\times T}(U)
        =\{(a,t)\in H_f^\circ\times T\mid a\in t^NU\}.
    \]
    The homomorphism $\psi:T_{f,N}\to S_f$ defined by
    $\psi(t)\coloneqq(t^N,t)$ induces an isomorphism on connected
    components, as we now show.

    Set $V_f\coloneqq U\cap\theta_f(U)$.
    Since $T$ is $\theta_f$-stable, we have a $\theta_f$-stable decomposition
    \[
        B\cap\theta_f(B)=T\ltimes V_f.
    \]
    Uniqueness of this decomposition
    therefore gives
    $B\cap H_f=T^{\theta_f}\ltimes V_f^{\theta_f}$.

    The group $V_f^{\theta_f}$ is connected. Indeed, choose a cocharacter
    $\nu$ of $T$ whose pairing with every root of $U$ is positive.
    If $\alpha$ is a root of $V_f$, both $\alpha$ and $\theta_f(\alpha)$
    are roots of $U$. Hence $\alpha$ has positive pairing with both
    $\nu$ and $\theta_f(\nu)$. Conjugation by the $\theta_f$-fixed cocharacter
    $\nu+\theta_f(\nu)$ contracts $V_f$ to $1$.
    This contraction preserves $V_f^{\theta_f}$, proving connectedness. In particular,
    $V_f^{\theta_f}\subseteq H_f^\circ$, and hence
    \[
        B\cap H_f^\circ
        =(T\cap H_f^\circ)\ltimes V_f^{\theta_f}.
    \]
    Projection onto the second coordinate therefore gives a split exact sequence
    \[
        1\longrightarrow V_f^{\theta_f}
        \xrightarrow{\ v\mapsto(v,1)\ }S_f
        \xrightarrow{\ (a,t)\mapsto t\ } T_{f,N}
        \longrightarrow1,
    \]
    with section $\psi$.
    Since $V_f^{\theta_f}$ is connected, $\psi$ induces
    $\pi_0(T_{f,N})\cong\pi_0(S_f)$.

    The group $H_f^\circ\times T$ is connected. Thus taking the stalk at $U$
    identifies irreducible equivariant local systems on $H_f^\circ B/U$
    with irreducible representations of $\pi_0(S_f)$.
    Since $T_{f,N}$ is a subgroup of the torus $T$, its component group is
    abelian, so these representations are characters.

    The local system associated with $(f,\mu)$ has deck character
    $\mu|_{T[N]}$, since $\psi(z)=(1,z)$ for $z\in T[N]$.
    By Definition~\ref{def:monodromic}, its torus monodromy is $\cL_{f,\mu}$.
    Definition~\ref{def:standard-basis-parameter} therefore identifies
    $\Gamma^{H^\circ}_{\cL'}$ with the pairs satisfying
    $\cL_{f,\mu}\cong\cL'$. Taking the union over $\cL'\in\fo$
    proves the bijection.
\end{proof}

\subsection{Cross action}
\label{subsec:cross-action}
\begin{definition}[cf. {\cite[\S7.5(c)--(c$'$)]{lusztig-vogan2014quasisplit}}]
    \label{def:cross-action}
    Let $w\in W$ and $(f,\mu)\in \Gamma^{H^\circ}_\fo$.
    For a positive system $P\subseteq\Phi(G,T)$, set
    \[
        \opr{Inv}_{f,P}(w)\coloneqq
        \{\beta\in P\mid\theta_f(\beta)=-\beta,\ w\beta\in-P\}.
    \]
    Let $\dot{w}\in N_G(T)$ be a lift of $w$.
    Let $w\times f$ be the unique element in $\mathscr{V}\cap H^\circ f\dot w^{-1}B$.
    Conjugation by $\dot w$ identifies $T_{f,N}$ with
    $T_{w\times f,N}$; write $(w\mu)(t)\coloneqq\mu(w^{-1}(t))$
    for the transported character.
    For $t\in T_{w\times f,N}$, set
    \[
        \begin{aligned}
            (w\times_{f,P}\mu)(t)
             & \coloneqq \mu(w^{-1}(t))
            \prod_{\beta\in\opr{Inv}_{w\times f,P}(w^{-1})}\beta(t^N), \\
            w\times_P(f,\mu)
             & \coloneqq (w\times f,w\times_{f,P}\mu).
        \end{aligned}
    \]
    Whenever the positive system is omitted from the notation,
    we mean $P=\Phi^+$, the positive system determined by $B$.
\end{definition}
\begin{lemma}
    \label{lem:cross-action-well-defined}
    For each positive system $P$, Definition~\ref{def:cross-action}
    defines an action $\times_P$ of $W$ on $\Gamma^{H^\circ}_\fo$.
    For $w\in W$ and $\xi=(f,\mu)$, one has
    \[
        \theta_{w\times f}=w\theta_fw^{-1}\quad\text{on }T,
        \qquad
        \cL_{w\times_P\xi}\cong w\cL_\xi.
    \]
    For $P=\Phi^+$, $s\in S$, and $t\in T_{s\times f,N}$,
    \[
        (s\times_f\mu)(t)=
        \begin{cases}
            \mu(s(t))\alpha_s(t^N), & \alpha_s\text{ is }\theta_f\text{-real}, \\
            \mu(s(t)),              & \text{otherwise}.
        \end{cases}
    \]
\end{lemma}
\begin{proof}
    Fix $P$ and write $\xi=(f,\mu)$.
    By the choice of $\mathscr V$ in
    \S\ref{subsec:classification-of-standard-objects},
    $w\times f$ lies in $H^\circ f\dot w^{-1}T$.
    Changing the lift of $w$ or the representative of this double coset
    does not change its orbit or the identification of subgroups of $T$.
    This proves the asserted conjugacy of involutions on $T$ and
    $T_{w\times f,N}=w(T_{f,N})$.
    It also gives $1\times f=f$ and
    $(vw)\times f=v\times(w\times f)$.

    If $\beta$ is $\theta_{w\times f}$-real and $t\in T_{w\times f,N}$,
    then $\beta(t^N)=\beta(t^N)^{-1}$.
    Hence each correction factor takes values in $\{\pm1\}$,
    is trivial on $T_{w\times f,N}^\circ$, and is trivial on $T[N]$.
    The corrected character therefore descends to the component group
    and has torus monodromy $w\cL_\xi$.
    Since $\fo$ is $W$-stable, $w\times_P\xi$ belongs to $\Gamma^{H^\circ}_\fo$.

    Let $v,w\in W$ and $t\in T_{(vw)\times f,N}$.
    For each real root pair $\{\beta,-\beta\}$ at $(vw)\times f$,
    the sign change relative to $P$ under $w^{-1}v^{-1}$ is the sum modulo two
    of the sign changes under $v^{-1}$ and $w^{-1}$.
    Since $\beta(t^N)=(-\beta)(t^N)$ and $\beta(t^N)^2=1$, this gives
    \[
        \begin{split}
            \prod_{\beta\in\opr{Inv}_{(vw)\times f,P}((vw)^{-1})}\beta(t^N)
            ={} & \prod_{\beta\in\opr{Inv}_{(vw)\times f,P}(v^{-1})}\beta(t^N) \\
                & \mathrel{}\cdot
            \prod_{\gamma\in\opr{Inv}_{w\times f,P}(w^{-1})}
            \gamma(v^{-1}(t)^N).
        \end{split}
    \]
    Together with $(vw)\mu=v(w\mu)$, this proves
    $(vw)\times_P\xi=v\times_P(w\times_P\xi)$.
    The identity has no inversions.
    Finally, for $P=\Phi^+$, a simple reflection has the single positive inversion
    $\alpha_s$, which is real at $s\times f$ exactly when it is real at $f$.
    This gives the displayed simple-reflection formula.
\end{proof}

\subsection{Root types and orbit geometry}
\label{subsec:root-types-and-orbit-geometry}
For $s\in S$, let $\alpha_s$ be its simple root.
\begin{definition}[cf. {\cite[\S7.5(d)--(g)]{lusztig-vogan2014quasisplit}}]
    \label{def:root-types}
    We say that the root $\alpha_s$ is
    \begin{enumerate}
        \item \emph{$\theta_f$-complex} if $\theta_f(\alpha_s)\ne\pm\alpha_s$;
              it is an ascent if $\theta_f(\alpha_s)\in\Phi^+$,
              and a descent if $\theta_f(\alpha_s)\in-\Phi^+$;
        \item \emph{$\theta_f$-real} if $\theta_f(\alpha_s)=-\alpha_s$;
        \item \emph{$\theta_f$-imaginary} if $\theta_f(\alpha_s)=\alpha_s$;
              it is compact or non-compact according as $d\theta_f$
              acts on $\mathfrak g_{\alpha_s}$ by $+1$ or $-1$.
    \end{enumerate}
\end{definition}

Let $P_s\supseteq B$ be the minimal parabolic subgroup associated with $s$.
The root type determines the $H^\circ$-orbits meeting $fP_s/B$.
Choose a root homomorphism $\phi_{s,f}:\SL_2\to P_s$ with the real-root
normalization of Lemma~\ref{lem:rank-one-involution-models} when $s$ is real.
Set $\dot{s}\coloneqq \phi_{s,f}\left(\begin{smallmatrix}
            0 & -1 \\
            1 & 0
        \end{smallmatrix}\right)\in N_G(T)$, which lifts $s$.
We have $P_s/B\cong \bP_1$ with affine charts
\[
    \begin{gathered}
        u_0(r) \coloneqq \phi_{s,f}\begin{pmatrix}
            1 & 0 \\
            r & 1
        \end{pmatrix},\qquad
        y_0(r) \coloneqq u_0(r)B,\\
        u_\infty(r) \coloneqq \phi_{s,f}\begin{pmatrix}
            1 & r \\
            0 & 1
        \end{pmatrix},\qquad
        y_\infty(r) \coloneqq u_\infty(r)\dot{s}B.
    \end{gathered}
\]
Identifying $P_s/B$ with the lines in $k^2$ gives
$y_0(r)=[1:r]$ and $y_\infty(r)=[r:1]$.
Let $\rho:P_s\to\Aut(P_s/B)\cong\PGL_2$ be the action map, and set
$Q\coloneqq\rho(H_f^\circ\cap P_s)$.
In these coordinates, $\rho(\phi_{s,f}(A))=[A]$ for $A\in\SL_2$.
Multiplication by $f$ identifies $Q$-orbits in $P_s/B$ with
$H^\circ$-orbits in $H^\circ fP_s/B$.

\begin{lemma}
    \label{lem:complex-ascent}
    If $s$ is $\theta_f$-complex ascent,
    then the $Q$-orbits in $P_s/B$ are $\{B\}$ and $(P_s\setminus B)/B$.
\end{lemma}
\begin{proof}
    The root set of $P_s$ is $\Phi^+\cup\{-\alpha_s\}$.
    Since $\theta_f(\alpha_s)$ is positive and distinct from $\alpha_s$,
    the root $-\alpha_s$ does not belong to $\theta_f(P_s)$.
    Hence $H_f^\circ\cap P_s\subseteq P_s\cap\theta_f(P_s)\subseteq B$,
    so $Q$ fixes $B$.

    Recall $V_f\coloneqq U\cap\theta_f(U)$.
    The action of $U$ on $(P_s\setminus B)/B$ is transitive,
    and factors through $\pi_s:U \to U/R_u(P_s)\cong \bG_a$.
    Therefore, it is enough to show
    that $d\pi_s$ is not trivial on $\opr{Lie}\left(V_f^{\theta_f}\right)$.
    We have
    \[
        X+d\theta_f(X)\in \opr{Lie}\left(V_f^{\theta_f}\right)
        ,\qquad 0\ne X \in \fg_{\alpha_s}.
    \]
    Since $\theta_f(\alpha_s)\ne\alpha_s$ is positive,
    $d\pi_s(d\theta_f(X))=0$, whereas $d\pi_s(X)\ne0$.
    Fixed subgroups are smooth in odd characteristic, so the image
    of $V_f^{\theta_f}$ in $\bG_a$ is all of $\bG_a$.
    This proves transitivity on the complement of $B$.
    Moreover, $\alpha_s+\theta_f(\alpha_s)\ne0$, so
    $\alpha_s$ restricts nontrivially to $(T^{\theta_f})^\circ$.
    Hence $\rho((T^{\theta_f})^\circ)=\rho(T)$ and $Q=\rho(B)$.
\end{proof}
\begin{lemma}
    \label{lem:complex-descent}
    If $s$ is $\theta_f$-complex descent,
    then the $Q$-orbits in $P_s/B$ are $\{\dot{s}B\}$ and $P_s/B\setminus\{\dot{s}B\}$.
\end{lemma}
\begin{proof}
    If $s$ is $\theta_f$-complex descent,
    then it is $\theta_{s\times f}$-complex ascent.
    So conjugating $Q$ by $\rho(\dot{s})$ proves the claim.
\end{proof}

For $\theta_f$-imaginary or real $s$
set $L_s\coloneqq\langle T,U_{\alpha_s},U_{-\alpha_s}\rangle$.
The group $L_s$ is a
$\theta_f$-stable common Levi subgroup of $P_s$ and $\theta_f(P_s)$.
We have
\[
    P_s\cap\theta_f(P_s)
    =L_s\ltimes\bigl(R_u(P_s)\cap\theta_f(R_u(P_s))\bigr).
\]
Both factors are $\theta_f$-stable, so the Levi component of every
element of $H_f^\circ\cap P_s$ belongs to $L_s^{\theta_f}$.
As $R_u(P_s)$ acts trivially on $P_s/B$, this gives
$Q\subseteq\rho(L_s^{\theta_f})$.
The restriction of $\rho$ identifies $L_s/Z(L_s)$ with $\PGL_2$.
Put
\[
    J_0\coloneqq\begin{pmatrix}1&0\\0&-1\end{pmatrix},
    \qquad J_1\coloneqq\begin{pmatrix}0&1\\1&0\end{pmatrix}.
\]

\begin{lemma}
    \label{lem:rank-one-involution-models}
    For imaginary or real $s$, the root homomorphism $\phi_{s,f}$ can be chosen so that
    $\theta_f\circ\phi_{s,f}=\phi_{s,f}\circ\widetilde\theta_f$, with
    \[
        \begin{array}{l|c|c}
            \text{root type}             & \widetilde\theta_f                                                & \SL_2^{\widetilde\theta_f} \\ \hline
            \text{compact imaginary}     & \mathrm{id}                                                       & \SL_2                      \\
            \text{non-compact imaginary} & \operatorname{Int}(J_0)
                                         & \{\operatorname{diag}(a,a^{-1}):a\in k^\times\}                                                \\
            \text{real}                  & \operatorname{Int}(J_1)
                                         & \left\{\left(\begin{smallmatrix}a&b\\b&a\end{smallmatrix}\right):
            a^2-b^2=1\right\}.
        \end{array}
    \]
    All three fixed subgroups are connected, and their images lie in $H_f^\circ$.
\end{lemma}
\begin{proof}
    The restriction of $\theta_f$ to $L_s^{\mathrm{der}}$ lifts uniquely
    through its simply connected covering $\phi_{s,f}:\SL_2\to L_s^{\mathrm{der}}$.
    For imaginary $s$, the lift fixes the diagonal torus pointwise.
    Its scalars on the two root subgroups are both $+1$ in the compact
    case and both $-1$ in the non-compact case. This gives the first two models.

    For real $s$, the lift exchanges the two root subgroups.
    It is therefore conjugation by
    $\left(\begin{smallmatrix}0&\kappa\\1&0\end{smallmatrix}\right)$
    for some $\kappa\in k^\times$.
    Choose $b\in k^\times$ with $b^2=\kappa$. Since
    \[
        \operatorname{diag}(b,1)^{-1}
        \begin{pmatrix}0&\kappa\\1&0\end{pmatrix}
        \operatorname{diag}(b,1)=bJ_1,
    \]
    replacing $\phi_{s,f}$ by
    $\phi_{s,f}\circ\operatorname{Int}(\operatorname{diag}(b,1))$
    gives the real model.
    The fixed subgroups follow by matrix multiplication.
\end{proof}

\begin{lemma}
    \label{lem:compact-imaginary-orbits}
    If $s$ is $\theta_f$-imaginary compact, then
    $Q$ acts transitively on $P_s/B$.
\end{lemma}
\begin{proof}
    By Lemma~\ref{lem:rank-one-involution-models},
    $\phi_{s,f}(\SL_2)\subseteq H_f^\circ\cap P_s$.
    Its projective image is $\PGL_2$, so $Q=\PGL_2$.
\end{proof}
\begin{lemma}
    \label{lem:noncompact-imaginary-orbits}
    If $s$ is $\theta_f$-imaginary non-compact, then
    $Q^\circ = \rho(T)$. Therefore, $P_s/B\setminus\{B,\dot{s}B\}$
    is a $Q$-orbit.
\end{lemma}
\begin{proof}
    The fixed diagonal torus in Lemma~\ref{lem:rank-one-involution-models}
    maps onto $\rho(T)$, since $k$ is algebraically closed.
    The Levi reduction therefore gives
    \[
        \rho(T)\subseteq Q\subseteq C_{\PGL_2}([J_0])=N_{\PGL_2}(\rho(T)).
    \]
    Indeed, an invertible matrix commutes projectively with $J_0$
    exactly when it preserves or exchanges the two eigenspaces of $J_0$.
    The normalizer has identity component $\rho(T)$, proving $Q^\circ=\rho(T)$.
    The torus is transitive on the complement of $\{B,\dot sB\}$,
    and the normalizer preserves that complement.
\end{proof}

For real or non-compact imaginary $s$, choose $\sqrt2\in k$ and put
\begin{equation}
    \label{eq:c-gs-def}
    c\coloneqq\frac1{\sqrt2}\begin{pmatrix}1&-1\\1&1\end{pmatrix},
    \qquad g_s\coloneqq\phi_{s,f}(c).
\end{equation}
Let $f_s$ be the unique element of $\mathscr V$ in $H^\circ fg_sB$.
Direct multiplication gives
\[
    c^{-1}J_0c=-J_1,\qquad c^{-1}J_1c=J_0,\qquad
    g_s^{-1}\theta_f(g_s)=\dot s^{-1}.
\]
Moreover, $g_s\dot s=\dot s g_s$, and the two columns of $c$ give
$g_sB=[1:1]$ and $g_s\dot sB=[-1:1]$.

\begin{lemma}
    \label{lem:real-root-orbits}
    If $s$ is $\theta_f$-real, then
    $Q^\circ = \rho(g_sTg_s^{-1})$.
    Therefore, $P_s/B\setminus\{g_sB,\dot{s}g_sB\}$
    is a $Q$-orbit.
\end{lemma}
\begin{proof}
    The fixed torus in Lemma~\ref{lem:rank-one-involution-models}
    maps onto $\rho(g_sTg_s^{-1})$.
    The Levi reduction and $c^{-1}J_1c=J_0$ give
    \[
        \rho(T)\subseteq\rho(g_s)^{-1}Q\rho(g_s)
        \subseteq C_{\PGL_2}([J_0]).
    \]
    Apply the normalizer and orbit argument in the proof of
    Lemma~\ref{lem:noncompact-imaginary-orbits} to
    $\rho(g_s)^{-1}Q\rho(g_s)$ and conjugate back by $\rho(g_s)$.
    The boundary points $B,\dot sB$ become $g_sB,g_s\dot sB$.
\end{proof}

\begin{definition}
    \label{def:type-i-ii}
    A $\theta_f$-real or non-compact imaginary root $\alpha_s$
    is of \emph{type I} if $Q$ fixes both boundary points of the open orbit individually
    and of \emph{type II} otherwise.
\end{definition}
The quotient $N_{\PGL_2}(Q^\circ)/Q^\circ$ has order two,
and its nontrivial coset exchanges the boundary points.
Thus $Q=Q^\circ$ in type I and $Q=N_{\PGL_2}(Q^\circ)$ in type II.

\begin{lemma}
    \label{lem:rank-one-cayley-geometry}
    Suppose $s$ is real or non-compact imaginary.
    The torus $g_sTg_s^{-1}$ is $\theta_f$-stable, and
    \[
        f_s\in H^\circ fg_sT,\qquad
        T\cap H_{f_s}^\circ=T\cap g_s^{-1}H_f^\circ g_s,\qquad
        \theta_{f_s}|_T=s\circ\theta_f|_T.
    \]
    The root $\alpha_s$ is non-compact imaginary for $\theta_{f_s}$
    if it is real for $\theta_f$, and real for $\theta_{f_s}$ otherwise.
    The type I or II of $\alpha_s$ is unchanged.
\end{lemma}
\begin{proof}
    The identity $g_s^{-1}\theta_f(g_s)=\dot s^{-1}$ shows that
    $g_sTg_s^{-1}$ is $\theta_f$-stable.
    The representative property in
    \S\ref{subsec:classification-of-standard-objects} gives
    $f_s=hfg_st_0$ for some $h\in H^\circ$ and $t_0\in T$.
    Hence
    \[
        H_{f_s}^\circ=t_0^{-1}g_s^{-1}H_f^\circ g_st_0,\qquad
        \theta_{f_s}
        =\operatorname{Int}(t_0^{-1}\dot s^{-1})\circ
        \theta_f\circ\operatorname{Int}(t_0).
    \]
    Intersecting the first equality with $T$ gives the stabilizer identity.
    Restricting the second equality to $T$ gives
    $\theta_{f_s}|_T=s\circ\theta_f|_T$.
    The identities $c^{-1}J_0c=-J_1$ and $c^{-1}J_1c=J_0$
    interchange the real and non-compact imaginary models.
    The additional conjugation by $t_0$ preserves these root types.

    Since $g_st_0\in P_s$, the two projective action groups satisfy
    \[
        \rho(H_{f_s}^\circ\cap P_s)
        =\rho(g_st_0)^{-1}Q\rho(g_st_0).
    \]
    Conjugation identifies their open orbits and boundary permutations.
    Definition~\ref{def:type-i-ii} therefore gives the same type for both.
\end{proof}

\begin{lemma}
    \label{lem:rank-one-boundary-permutation-sign}
    Suppose $s$ is $\theta_f$-real.
    Let $T_{f,N}$ act on $P_s/B$ by $t\cdot_N xB\coloneqq t^NxB$.
    For the boundary points $g_sB$ and $g_s\dot sB$, one has
    \[
        \delta(t)\coloneqq\alpha_s(t^N)=
        \begin{cases}
            1,  & \text{if $t$ fixes both boundary points},    \\
            -1, & \text{if $t$ exchanges the boundary points}.
        \end{cases}
    \]
    Some $t\in T_{f,N}$ exchanges the boundary points exactly
    when $s$ is of type II.
\end{lemma}
\begin{proof}
    Since $\theta_f(\alpha_s)=-\alpha_s$, the root $\alpha_s$
    does not occur in $V_f=U\cap\theta_f(U)$.
    Thus $V_f\subseteq R_u(P_s)$ acts trivially on $P_s/B$.
    The stabilizer decomposition in the proof of
    Lemma~\ref{lem:combi-parameter-rel-standard-obj} gives
    \[
        \operatorname{Stab}_Q(B)
        =\rho(H_f^\circ\cap B)=\rho(T\cap H_f^\circ).
    \]
    The torus $\rho(g_sTg_s^{-1})$ acts freely away from its two fixed points.
    The element $[J_0]$ fixes $B$ and exchanges $[1:1]$ and $[-1:1]$.
    The real model therefore gives
    \[
        \operatorname{Stab}_Q(B)=
        \begin{cases}
            \{1\},       & \text{type I},  \\
            \{1,[J_0]\}, & \text{type II}.
        \end{cases}
    \]
    By definition of $T_{f,N}$, surjectivity of $[N]:T\to T$
    gives $[N](T_{f,N})=T\cap H_f^\circ$.
    In the chart $y_\infty$, the action is
    $r\mapsto\alpha_s(t^N)r$, and the boundary points are $r=\pm1$.
    The two stabilizer elements act by $r\mapsto r$ and $r\mapsto-r$,
    proving both assertions.
\end{proof}

\begin{corollary}
    \label{cor:rank-one-fixed-points}
    The $H^\circ$-orbits in $H^\circ fP_s/B$ have the following indices in $\mathscr{V}$.
    \[
        \begin{array}{l|c|c|c}
            \text{root type}                     & Q\text{-orbits index} & \multicolumn{2}{c}{\text{cross-action relations}}                      \\ \hline
            \text{complex ascent}                & f,s\times f           & f\ne s\times f                                    & -                  \\
            \text{complex descent}               & f,s\times f           & f\ne s\times f                                    & -                  \\
            \text{imaginary compact}             & f                     & f=s\times f                                       & -                  \\
            \text{imaginary non-compact type I}  & f,s\times f,f_s       & f\ne s\times f                                    & f_s=s\times f_s    \\
            \text{imaginary non-compact type II} & f,f_s                 & f=s\times f                                       & f_s=s\times f_s    \\
            \text{real type I}                   & f,f_s,s\times f_s     & f=s\times f                                       & f_s\ne s\times f_s \\
            \text{real type II}                  & f,f_s                 & f=s\times f                                       & f_s=s\times f_s    \\
        \end{array}
    \]
\end{corollary}
\begin{proof}
    Apply Lemmas~\ref{lem:complex-ascent}--\ref{lem:real-root-orbits}.
    For real roots, $g_s\dot s=\dot s g_s$ identifies the second
    boundary orbit with $s\times f_s$.
    Lemma~\ref{lem:rank-one-cayley-geometry} identifies the root type at $f_s$.
    Applying the same orbit descriptions at $f_s$ gives the remaining relations.
\end{proof}

\subsection{Cayley transform and parity}
\label{subsec:cayley-transform-and-parity}

Choose a primitive $2N$-th root $\zeta_{2N}\in k^\times$ and put
$\zeta_N\coloneqq\zeta_{2N}^2$.
For a real root, the connected fixed torus in
Lemma~\ref{lem:rank-one-involution-models} contains $-I$.
Hence $\alpha_s^\vee(-1)\in T\cap H_f^\circ$ and
$\alpha_s^\vee(\zeta_{2N})\in T_{f,N}$.

\begin{definition}[cf. {\cite[\S7.5(g)]{lusztig-vogan2014quasisplit}}]
    \label{def:root-parity}
    For a real pair $(s,\xi)$, with $\xi=(f,\mu)$, put
    \[
        \eta_{s,\xi}\coloneqq\mu(\alpha_s^\vee(\zeta_{2N})).
    \]
    The pair is \emph{parity} if $\eta_{s,\xi}=1$ and
    \emph{non-parity} otherwise.
\end{definition}
Replacing $\zeta_{2N}$ by another primitive root raises
$\eta_{s,\xi}$ to a power coprime to $2N$.
Thus parity is independent of that choice.
By Lemma~\ref{lem:connected-center-character-weyl-group},
\[
    \eta_{s,\xi}^2=1
    \quad\Longleftrightarrow\quad
    (\alpha_s^\vee)^*\cL_\xi\cong\overline{\bQ}_\ell
    \quad\Longleftrightarrow\quad s\cL_\xi\cong\cL_\xi.
\]
In particular, a real pair with $s\cL_\xi\not\cong\cL_\xi$
is non-parity.

\begin{lemma}
    \label{lem:reflection-on-parameter-characters}
    If $s$ is imaginary or real parity, then $s(T_{f,N})=T_{f,N}$ and
    $\mu\circ s=\mu$.
    If a real pair has $\eta_{s,\xi}=-1$, then
    $s\times\xi=\xi$.
\end{lemma}
\begin{proof}
    Use $s(t)=t\alpha_s^\vee(\alpha_s(t)^{-1})$.
    For imaginary $s$, the torus $\alpha_s^\vee(\bG_m)$ lies in
    $T_{f,N}^\circ$, so its character under $\mu$ is trivial.
    For real $s$, one has $s\times f=f$ and
    $\alpha_s(t)^{2N}=1$ for $t\in T_{f,N}$.
    Write $\alpha_s(t)=\zeta_{2N}^a$. Then
    \[
        \mu(s(t))=\mu(t)\eta_{s,\xi}^{-a},\qquad
        \alpha_s(t^N)=(-1)^a.
    \]
    Parity gives $\mu(s(t))=\mu(t)$.
    If $\eta_{s,\xi}=-1$, the correction factor in
    Lemma~\ref{lem:cross-action-well-defined} cancels
    $\eta_{s,\xi}^{-a}$.
\end{proof}

For $s\in S$ and $f\in\mathscr V$, put
\[
    T_{s,f,N}\coloneqq\{t\in T_{f,N}\mid\alpha_s(t^N)=1\}.
\]
\begin{lemma}
    \label{lem:cayley-common-stabilizer}
    Suppose $s$ is $\theta_f$-non-compact imaginary. Then
    \[
        T_{s,f,N}=T_{s,f_s,N}.
    \]
    The subgroup $T_{s,f,N}$ contains $T[N]$.
    Moreover, $T_{f_s,N}^\circ\subseteq T_{f,N}^\circ$.
    The index $[T_{f_s,N}:T_{s,f,N}]$
    is one in type I and two in type II.
    For the $T_{f_s,N}$-action
    $t\cdot_{s,N}xB\coloneqq g_st^Ng_s^{-1}xB$
    on $\{B,\dot sB\}$, the stabilizer of either point is
    $T_{s,f,N}$.
\end{lemma}
\begin{proof}
    If $t\in T$ satisfies $\alpha_s(t^N)=1$, then $t^N$ centralizes
    $\phi_{s,f}(\SL_2)$ and hence $g_s$.
    Lemma~\ref{lem:rank-one-cayley-geometry} therefore gives
    $t^N\in H_{f_s}^\circ$ if and only if $t^N\in H_f^\circ$.
    This proves the subgroup equality and the containment of $T[N]$.

    By Lemmas~\ref{lem:rank-one-cayley-geometry} and
    \ref{lem:rank-one-boundary-permutation-sign},
    $\alpha_s(t^N)$ on $T_{f_s,N}$ has image $\{1\}$ in type I
    and $\{1,-1\}$ in type II.
    Its kernel is $T_{s,f,N}$, proving the index assertion.
    This kernel contains $T_{f_s,N}^\circ$, which is therefore
    a connected subgroup of $T_{f,N}$.
    Finally, $\rho(t^N)$ is $1$ or $[J_0]$ according to this sign.
    Conjugation by $\rho(g_s)=[c]$ sends $[J_0]$ to $[J_1]$,
    which exchanges $B$ and $\dot sB$. This proves the stabilizer assertion.
\end{proof}

\begin{definition}[cf. {\cite[\S7.5(h)--(k)]{lusztig-vogan2014quasisplit}}]
    \label{def:cayley-transform}
    Let $\xi=(f,\mu)$ and suppose $s$ is $\theta_f$-non-compact imaginary.
    The \emph{Cayley transform} of $\xi$ is
    \[
        c_s(\xi)\coloneqq
        \left\{(f_s,\mu'):
        \mu'\in\widehat{\pi_0(T_{f_s,N})},\quad
        \mu'|_{T_{s,f,N}}
        =\mu|_{T_{s,f,N}}\right\}.
    \]
    For a real parameter $\zeta$, its \emph{inverse Cayley transform} is
    \[
        c^s(\zeta)\coloneqq
        \{\xi:(s,\xi)\text{ is non-compact imaginary},\
        \zeta\in c_s(\xi)\}.
    \]
\end{definition}
For every $\zeta=(f_s,\mu')\in c_s(\xi)$, the restriction condition gives
$\mu'|_{T[N]}=\mu|_{T[N]}$, so
$\cL_\zeta\cong\cL_\xi$.

\begin{lemma}
    \label{lem:cayley-transform-computation}
    Cayley transforms are real parity parameters of the same type.
    Their cardinalities are
    \[
        \begin{aligned}
            |c_s(\xi)| & =
            \begin{cases}
                1, & \text{non-compact imaginary I},  \\
                2, & \text{non-compact imaginary II},
            \end{cases} \\
            |c^s(\xi)| & =
            \begin{cases}
                2, & \text{real parity I},   \\
                1, & \text{real parity II},  \\
                0, & \text{real non-parity}.
            \end{cases}
        \end{aligned}
    \]
    Every two-element set consists of distinct parameters exchanged
    by cross action. Cross action fixes non-compact imaginary type II
    and real parity type I parameters.
\end{lemma}
\begin{proof}
    By Lemma~\ref{lem:rank-one-cayley-geometry}, $\alpha_s$ is
    $\theta_{f_s}$-real and has the same type as the
    $\theta_f$-non-compact imaginary root $\alpha_s$.
    Lemma~\ref{lem:cayley-common-stabilizer} gives
    $T_{f_s,N}^\circ\subseteq T_{f,N}^\circ$.
    Hence $\mu|_{T_{s,f,N}}$ defines a character of a subgroup of
    $\pi_0(T_{f_s,N})$. This finite group is abelian, so the character
    extends. Its extensions form a torsor under the character group
    of $T_{f_s,N}/T_{s,f,N}$.
    By Lemma~\ref{lem:cayley-common-stabilizer}, the quotient
    has order one in type I and two in type II.
    This proves the asserted cardinalities of $c_s(\xi)$.
    Since $\alpha_s^\vee(\zeta_{2N})\in T_{s,f,N}$ and
    $\alpha_s^\vee(\bG_m)\subseteq T_{f,N}^\circ$, all extensions
    satisfy parity.

    Conversely, fix a real parameter $(f_s,\mu')$ and a boundary
    orbit represented by $f$, so $\alpha_s$ is $\theta_f$-non-compact imaginary.
    By Lemma~\ref{lem:cayley-common-stabilizer},
    $T_{s,f,N}$ is the kernel of
    $t\mapsto\alpha_s(t)^N$ on $T_{f,N}$.
    Since $\alpha_s^\vee(\bG_m)\subseteq T_{f,N}^\circ$ and
    $\alpha_s(\alpha_s^\vee(z))^N=z^{2N}$, every component of
    $T_{f,N}$ meets $T_{s,f,N}$.
    Thus inclusion induces a surjection
    \[
        \pi_0(T_{s,f,N})\longrightarrow\pi_0(T_{f,N}).
    \]

    Its kernel consists of components represented by
    $t\in T_{s,f,N}\cap T_{f,N}^\circ$.
    The identity component of the kernel of
    $\alpha_s|_{T_{f,N}^\circ}$ has codimension one.
    Since $\alpha_s\circ\alpha_s^\vee(z)=z^2$, multiplication from
    this connected kernel and $\alpha_s^\vee(\bG_m)$ onto
    $T_{f,N}^\circ$ is surjective.
    Write $t=\alpha_s^\vee(z)t_0$, where
    $t_0\in(\ker(\alpha_s|_{T_{f,N}^\circ}))^\circ$.
    Then $z^{2N}=1$ and $t_0\in(T_{s,f,N})^\circ$.
    Hence the kernel on component groups is generated by the
    component of $\alpha_s^\vee(\zeta_{2N})$.

    Therefore $\mu'|_{T_{s,f,N}}$ descends uniquely to
    a character of $\pi_0(T_{f,N})$ exactly when
    $\mu'(\alpha_s^\vee(\zeta_{2N}))=1$, namely when $(f_s,\mu')$ is parity.
    Type I has two boundary orbits, and type II has one.
    This proves the inverse-transform assertions.

    Lemma~\ref{lem:reflection-on-parameter-characters} shows that
    cross action transports the same boundary character between the
    two type I boundary orbits and fixes a type II boundary parameter.
    On a real parity parameter, cross action multiplies its character
    by the boundary sign $\delta(t)=\alpha_s(t^N)$ of
    Lemma~\ref{lem:rank-one-boundary-permutation-sign}.
    It therefore fixes type I and exchanges the two extensions in type II.
\end{proof}

\section{Formulas for the simple Hecke operators}
\label{sec:formulas-for-the-simple-hecke-operators}
We compute the coefficient $a_{s,\xi,\xi'}$ of $m_{\xi'}$ in $\mathbf T_s m_\xi$, where
$s\in S$, $\xi=(f,\mu)\in\Gamma^{H^\circ}_\fo$ is the input parameter,
and $\xi'=(f',\mu')\in\Gamma^{H^\circ}_\fo$ is the output parameter.
By \eqref{eq:convolution-formula-for-hecke-action}, the coefficients are determined
by the weight-graded cohomology
of $(K_{s,\xi})_{f'U}$, viewed as a $T_{f',N}$-representation.
We identify each stalk with the compactly supported cohomology of
a punctured orbit in $P_s/B$. The orbit geometry,
local monodromy, and $T_{f',N}$-character determine every coefficient.
\subsection{Reduction to the support}
\label{subsec:reduction-to-support}
By \eqref{eq:convolution-formula-for-hecke-action},
\[
    \mathbf T_sm_\xi = \chi_u(K_{s,\xi}),\qquad
    K_{s,\xi} \coloneqq \overline p_{1!}q_{2\flat}
    \bigl(p_2^*(j_{\xi!}\mathscr E_\xi)
    \otimes j_{s!}\mathscr E_{s,\cL_\xi}\bigr).
\]
The stalk $(K_{s,\xi})_{xU}$ can be nonzero only if there is
$yB\in H^\circ fB/B$ with $\opr{pos}(xB,yB)=s$.
Corollary~\ref{cor:rank-one-fixed-points} lists the possible output orbits.
Using the first definition applicable for $f'$, set
\begin{equation}
    \label{eq:rank-one-output-representatives}
    \begin{gathered}
        g_{f'}\coloneqq
        \begin{cases}
            1, & f'=f,           \\
            \left(\begin{smallmatrix}0&-1\\1&0\end{smallmatrix}\right),
               & f'=s\times f,   \\
            c, & f'=f_s,         \\
            c\left(\begin{smallmatrix}0&-1\\1&0\end{smallmatrix}\right),
               & f'=s\times f_s.
        \end{cases},
        \quad
        p_{f'}\coloneqq\phi_{s,f}(g_{f'})
        ,\quad
        c\coloneqq\frac1{\sqrt2}\begin{pmatrix}1&-1\\1&1\end{pmatrix} \text{ as in \eqref{eq:c-gs-def}}.
    \end{gathered}
\end{equation}
By definition, $f'\in H^\circ fp_{f'}T$ and
$T\cap H_{f'}^\circ=T\cap H_{fp_{f'}}^\circ$.
$H^\circ$-equivariance identifies the geometric stalks at $f'U$ and
$fp_{f'}U$ with their $T_{f',N}$-actions.
Enlarge the auxiliary finite field so that a chosen transporting
element of $H^\circ\times T$ is rational. The stalk identification
then respects Frobenius.
Set
\[
    z_{f'}\coloneqq p_{f'}B,\qquad D_{f'}\coloneqq Q\cdot B\setminus\{z_{f'}\}.
\]
\begin{corollary}
    \label{cor:rank-one-punctured-orbits}
    The domains are given, up to isomorphism, by the following table.
    A dash denotes an absent output orbit.
    \[
        \begin{array}{l|c|c|c}
            \text{root type}                & D_f                    & D_{s\times f}\ (s\times f\ne f) & \text{other domains}              \\\hline
            \text{complex ascent}           & \emptyset              & \{B\}                           & -                                 \\
            \text{complex descent}          & \bG_m                  & \bA_1                           & -                                 \\
            \text{compact imaginary}        & \bA_1                  & -                               & -                                 \\
            \text{non-compact imaginary I}  & \emptyset              & \{B\}                           & D_{f_s}=\{B\}                     \\
            \text{non-compact imaginary II} & \{\dot sB\}            & -                               & D_{f_s}=\{B,\dot sB\}             \\
            \text{real I}                   & \bA_1\setminus\{1,-1\} & -                               & D_{f_s},D_{s\times f_s}\cong\bG_m \\
            \text{real II}                  & \bA_1\setminus\{1,-1\} & -                               & D_{f_s}\cong\bG_m
        \end{array}
    \]
\end{corollary}
\begin{proof}
    Remove $z_{f'}$ from the input orbit $Q\cdot B$ in
    Corollary~\ref{cor:rank-one-fixed-points}.
    For real inputs, the affine coordinate is $y_\infty$.
\end{proof}

For $t\in T_{f',N}$, set
\[
    h_t\coloneqq fp_{f'}t^Np_{f'}^{-1}f^{-1}\in H^\circ.
\]
The pair $(h_t,t)$ fixes $fp_{f'}U$ for the action
$(a,t)\cdot_N gU=agt^{-N}U$.
It acts on $\widetilde{\cB}\times\cB$ by
\[
    (g_1U,g_2B)\longmapsto(h_tg_1t^{-N}U,h_tg_2B).
\]
Identify the fiber of $\overline p_1$ over $fp_{f'}U$ with $\cB$ via
\[
    i_{f'}:\cB\longrightarrow\widetilde{\cB}\times\cB,
    \qquad xB\longmapsto(fp_{f'}U,fxB).
\]
The map $i_{f'}$ is equivariant for
\[
    t\cdot_{f',N} xB\coloneqq p_{f'}t^Np_{f'}^{-1}xB.
\]
Since $p_{f'}t^Np_{f'}^{-1}\in H_f^\circ\cap P_s$, this action
preserves $D_{f'}$.
The following square is Cartesian and $T_{f',N}$-equivariant:
\[
    \begin{tikzcd}
        \cB
        \ar[r,"i_{f'}"]
        \ar[d]
        \ar[dr,phantom,"\lrcorner",very near start]
        & \widetilde{\cB}\times \cB
        \ar[d,"\overline p_1"]
        \\
        \{fp_{f'}U\}
        \ar[r]
        & \widetilde{\cB}
    \end{tikzcd}
\]
Set
\[
    \mathscr F_{f'}\coloneqq
    (i_{f'}|_{D_{f'}})^*q_{2\flat}
    \bigl(p_2^*(j_{\xi!}\mathscr E_\xi)
    \otimes j_{s!}\mathscr E_{s,\cL_\xi}\bigr).
\]
The sheaf $\mathscr F_{f'}$ is a rank-one local system pure of weight zero,
since both factors defining it have these properties on the inverse image of $D_{f'}$.

Proper base change, followed by restriction to the support, gives
\[
    (K_{s,\xi})_{f'U}\cong(K_{s,\xi})_{fp_{f'}U}\cong
    R\Gamma_c(D_{f'},\mathscr F_{f'}).
\]
These isomorphisms are $T_{f',N}$-equivariant and compatible with Frobenius.
The action on cohomology factors through $\pi_0(T_{f',N})$.
Give each output character $\mu'$ the weight-zero Weil structure
of its local-system stalk. By Lemma~\ref{lem:combi-parameter-rel-standard-obj},
the coefficient for $\xi'=(f',\mu')$ is
\[
    a_{s,\xi,\xi'}=
    \sum_{i,j}(-1)^i(-u^{1/2})^j
    \dim\Hom_{\pi_0(T_{f',N})}\!\left(
    \mu',\opr{Gr}^{\mathrm{Wt}}_jH_c^i(D_{f'},\mathscr F_{f'})\right).
\]

\begin{lemma}
    \label{lem:rank-one-coefficient-restriction}
    On the domain where $y_{f'}(r)\coloneqq p_{f'}y_\infty(r)\in D_{f'}$, put
    \[
        \tilde y_{f'}(r)\coloneqq
        p_{f'}u_\infty(r)\dot sU.
    \]
    \begin{enumerate}[label=(\roman*),ref=(\roman*)]
        \item \label{itm:rank-one-geometric-coefficients}
              There is a noncanonical geometric isomorphism
              \[
                  (y_{f'})^*\mathscr F_{f'}
                  \cong (f\tilde y_{f'})^*\mathscr E_\xi.
              \]
              The local system $\mathscr F_{f'}^{\otimes 2N}$ is geometrically
              trivial. In particular, $\mathscr F_{f'}$ is tame.
        \item \label{itm:rank-one-stabilizer-character}
              For $xB\in D_{f'}$, the stabilizer
              $\opr{Stab}_{T_{f',N}}(xB)$ acts on
              $(\mathscr F_{f'})_{xB}$ by the character
              \[
                  t\longmapsto\mu(s(t)).
              \]
    \end{enumerate}
\end{lemma}
\begin{proof}
    For $xB\in D_{f'}$, there is a unique $r\in k$ such that
    $xB=y_{f'}(r)$. Choose the representative
    $x\coloneqq p_{f'}u_\infty(r)\dot s$.
    Descent along $q_2$ identifies
    \[
        (\mathscr F_{f'})_{xB}\cong
        (\mathscr E_\xi)_{fxU}\otimes
        (\mathscr E_{s,\cL_\xi})_{(fp_{f'}U,fxU)}.
    \]
    Diagonal multiplication by $(fp_{f'}u_\infty(r))^{-1}$
    carries $(fp_{f'}U,fxU)$ to $(U,\dot sU)$.

    \emph{\ref{itm:rank-one-geometric-coefficients}.}
    As $\mathscr E_{s,\cL_\xi}$ is $G$-equivariant, its pullback
    by $\tilde{y}_{f',f}(r) \coloneqq \left(fp_{f'}U,f\tilde y_{f'}(r)\right)$
    is geometrically constant.
    Trivializing $\tilde{y}_{f',f}^*\mathscr E_{s,\cL_\xi}$ gives
    the isomorphism in \ref{itm:rank-one-geometric-coefficients}.

    The image of $T\to T$, $x\mapsto x\theta_f(x)$, is connected
    and contained in $T_{f,N}^\circ$.
    For $t\in T_{f,N}$, it contains
    $t^{2N}=t^N\theta_f(t^N)$. Hence $\mu^{2N}=1$.
    Lemma~\ref{lem:combi-parameter-rel-standard-obj} gives
    $\mathscr E_\xi^{\otimes2N}\cong\overline{\bQ}_\ell$
    geometrically. Since $y_{f'}$ identifies its domain with $D_{f'}$,
    the isomorphism in \ref{itm:rank-one-geometric-coefficients} gives
    geometric triviality of $\mathscr F_{f'}^{\otimes2N}$.
    Since $2N$ is prime to the
    characteristic, wild inertia acts trivially.

    \emph{\ref{itm:rank-one-stabilizer-character}.}
    Let $t\in\opr{Stab}_{T_{f',N}}(xB)$.
    The identity
    \[
        t^Nu_\infty(r)t^{-N}
        =u_\infty(\alpha_s(t)^Nr)
    \]
    gives $\alpha_s(t)^Nr=r$. Hence $t^N$ commutes with
    $u_\infty(r)$, and
    \[
        h_tfp_{f'}t^{-N}U=fp_{f'}U,\qquad
        h_tfxs(t)^{-N}U=fxU.
    \]
    $H^\circ\times T$-equivariance of descent expresses the action of $t$ as the tensor
    product of the actions of $(h_t,s(t))$ on $(\mathscr E_\xi)_{fxU}$ and
    $(h_t,t,s(t))$ on $(\mathscr E_{s,\cL_\xi})_{(fp_{f'}U,fxU)}$.

    Since $xB\in Q\cdot B$, some element $(h_0,t_0)\in H^\circ\times T$
    carries $fxU$ to $fU$ for the $N$-th power action.
    Conjugation by $(h_0,t_0)$ preserves the torus coordinate $s(t)$.
    The stabilizer description in the proof of
    Lemma~\ref{lem:combi-parameter-rel-standard-obj} therefore gives
    $s(t)\in T_{f,N}$ and the input character $\mu(s(t))$.

    Diagonal multiplication by $(fp_{f'}u_\infty(r))^{-1}$
    conjugates $(h_t,t,s(t))$ to $(t^N,t,s(t))$.
    The latter belongs to the connected subgroup
    \[
        \{(a^N,a,s(a))\mid a\in T\}
        \subseteq\opr{Stab}_{G\times T\times T}(U,\dot sU)
    \]
    for the $N$-th power action. This subgroup acts trivially on
    $(\mathscr E_{s,\cL_\xi})_{(U,\dot sU)}$.
    Thus $(h_t,t,s(t))$ acts trivially on
    $(\mathscr E_{s,\cL_\xi})_{(fp_{f'}U,fxU)}$, and the character
    on $(\mathscr F_{f'})_{xB}$ is $\mu(s(t))$.
\end{proof}

\subsection{Complex roots}
\label{subsec:complex-roots}
\begin{lemma}
    \label{lem:complex-descent-monodromy}
    For a complex descent, $\mathscr F_f$ on
    $D_f=\bP_1\setminus\{B,\dot sB\}$ is geometrically constant
    if and only if $s\cL_\xi\cong\cL_\xi$.
    In that case $T_{f,N}$ acts on $H^0(D_f,\mathscr F_f)$ by $\mu$.
\end{lemma}
\begin{proof}
    For $z\ne0$, direct matrix multiplication in $\SL_2$, followed by
    application of $\phi_{s,f}$, gives
    \[
        u_\infty(z^{-1})\dot s
        =u_0(z)\alpha_s^\vee(z^{-1})u_\infty(-z).
    \]
    Since $u_\infty(-z)\in U$, we obtain
    \[
        fu_\infty(z^{-1})\dot sU
        =fu_0(z)\alpha_s^\vee(z^{-1})U.
    \]
    The map $z\mapsto fu_0(z)U$ takes values in
    $\widetilde{\cO}_\xi$ for all $z\in\bA_1$, since the input
    orbit omits only $\dot sB$ in this fiber.
    Therefore $(z\mapsto fu_0(z)U)^*\mathscr E_\xi$ extends across $z=0$.
    By Lemma~\ref{lem:rank-one-coefficient-restriction}\ref{itm:rank-one-geometric-coefficients}, the monodromy
    of $\mathscr F_f$ at $B$ is therefore that of
    $(\alpha_s^\vee)^*\cL_\xi$, up to inversion.
    The two monodromies on $\bG_m$ are inverse.
    Lemma~\ref{lem:connected-center-character-weyl-group} now proves
    the constancy criterion.

    Suppose now that $s\cL_\xi\cong\cL_\xi$.
    The character $\alpha_s$ restricts nontrivially to
    $(T^{\theta_f})^\circ$, because $\theta_f(\alpha_s)\ne-\alpha_s$.
    Thus $T_{f,N}^\circ$ is transitive on $D_f$.
    A point stabilizer consists of the $t\in T_{f,N}$ with
    $\alpha_s(t)^N=1$.
    Trivial coroot monodromy implies
    $\mu(\alpha_s^\vee(\alpha_s(t)))=1$, so
    Lemma~\ref{lem:rank-one-coefficient-restriction}\ref{itm:rank-one-stabilizer-character} gives the stalk character
    $\mu(s(t))=\mu(t)$.
    Lemma~\ref{lem:global-section-character-from-stabilizer}
    therefore identifies the character on $H^0(D_f,\mathscr F_f)$ with $\mu$.
\end{proof}

\begin{proposition}
    \label{prop:complex-root-action}
    For a complex root,
    \[
        \mathbf T_sm_\xi=
        \begin{cases}
            m_{s\times\xi}, & \text{ascent},                               \\
            (u-1)m_\xi+u m_{s\times\xi},
                            & \text{descent and }s\cL_\xi\cong\cL_\xi,     \\
            u m_{s\times\xi},
                            & \text{descent and }s\cL_\xi\not\cong\cL_\xi.
        \end{cases}
    \]
\end{proposition}
\begin{proof}
    In an ascent, the only nonempty domain is
    $D_{s\times f}=\{B\}$.
    Lemma~\ref{lem:rank-one-coefficient-restriction}\ref{itm:rank-one-stabilizer-character} identifies its weight-zero
    stalk character with $\mu\circ s$, giving $m_{s\times\xi}$.
    In a descent, $D_{s\times f}\cong\bA_1$ contains the fixed
    point $B$. Its tame coefficient system is constant.
    Lemma~\ref{lem:affine-line-cohomology} gives a single degree-two,
    weight-two line with character $\mu\circ s$, contributing
    $u m_{s\times\xi}$.

    The remaining domain is $D_f\cong\bG_m$.
    Lemma~\ref{lem:complex-descent-monodromy} determines whether
    its coefficients are constant and identifies their character.
    Lemma~\ref{lem:multiplicative-group-cohomology} gives $(u-1)m_\xi$ in the
    constant case and zero otherwise.
\end{proof}

\subsection{Imaginary roots}
\label{subsec:imaginary-roots}
\begin{proposition}
    \label{prop:imaginary-root-action}
    For an imaginary root, $s\cL_\xi\cong\cL_\xi$ and
    \[
        \mathbf T_sm_\xi=
        \begin{cases}
            u m_\xi, & \text{compact},     \\
            m_{s\times\xi}+\displaystyle\sum_{\zeta\in c_s(\xi)}m_\zeta,
                     & \text{non-compact}.
        \end{cases}
    \]
\end{proposition}
\begin{proof}
    Lemma~\ref{lem:reflection-on-parameter-characters} gives $\mu\circ s=\mu$.
    Restriction to $T[N]$ gives $s\cL_\xi\cong\cL_\xi$ by
    Lemma~\ref{lem:combi-parameter-rel-standard-obj}.
    In the compact case, $D_f\cong\bA_1$ and $\dot sB$ is fixed
    by $T_{f,N}$. Lemma~\ref{lem:rank-one-coefficient-restriction}\ref{itm:rank-one-stabilizer-character} gives the
    character $\mu\circ s=\mu$ at that point.
    Lemma~\ref{lem:affine-line-cohomology} therefore gives $u m_\xi$.

    In non-compact type I, $D_f=\emptyset$ and
    $D_{s\times f}=\{B\}$. The latter contributes $m_{s\times\xi}$.
    In type II, the input $Q$-orbit is $\{B,\dot sB\}$.
    Removing $B$ leaves $D_f=\{\dot sB\}$, with character $\mu$.
    Its contribution is $m_\xi$.

    It remains to compute the coefficient over $f_s$.
    The domain $D_{f_s}$ has one point in type I and two points
    in type II. By Lemma~\ref{lem:cayley-common-stabilizer},
    $T_{f_s,N}$ acts transitively on this domain with point stabilizer
    $T_{s,f,N}$.
    Lemmas~\ref{lem:rank-one-coefficient-restriction}\ref{itm:rank-one-stabilizer-character} and
    \ref{lem:reflection-on-parameter-characters} give the stalk
    character $\mu|_{T_{s,f,N}}$.
    Thus its degree-zero cohomology is the induced representation
    \[
        \operatorname{Ind}_{T_{s,f,N}}^{T_{f_s,N}}
        \bigl(\mu|_{T_{s,f,N}}\bigr)
        \cong\bigoplus_{(f_s,\mu')\in c_s(\xi)}\mu'.
    \]
    Every summand has weight zero and multiplicity one.
    In type II, $s\times\xi=\xi$ by
    Lemma~\ref{lem:cayley-transform-computation}.
    Thus both non-compact cases give the stated sum.
\end{proof}

\subsection{Real roots}
\label{subsec:real-roots}
\begin{lemma}
    \label{lem:real-root-monodromy}
    Suppose $s$ is real and $D_{f'}\ne\emptyset$.
    In the output chart $p_{f'}y_\infty$, the local monodromies
    of $\mathscr F_{f'}$ are, up to simultaneous inversion,
    \[
        \begin{array}{c|c|c}
            \text{output} & \text{finite punctures}     & z_{f'}=p_{f'}B    \\\hline
            f'=f          & (\eta_{s,\xi},\eta_{s,\xi}) & \eta_{s,\xi}^{-2} \\
            f'\ne f       & \eta_{s,\xi}                & \eta_{s,\xi}^{-1}
        \end{array}
    \]
    In particular, $\mathscr F_{f'}$ is geometrically constant
    exactly when $(s,\xi)$ is parity.
\end{lemma}
\begin{proof}
    Identify $\SL_2/U_+$ with $k^2\setminus\{0\}$ by the first column.
    The real fixed torus in Lemma~\ref{lem:rank-one-involution-models}
    preserves $q(x,y)\coloneqq x^2-y^2$.
    For $g\in\SL_2$, put
    \[
        h_g(r)\coloneqq q(g(r,1)^{\mathsf t}).
    \]
    In the chart $y_{f'}$, the domain $D_{f'}$ is given by $h_{g_{f'}}(r)\ne0$.
    Here $g_{f'}$ is the matrix defined in \eqref{eq:rank-one-output-representatives}.
    Define
    \[
        \widetilde D_{f'}\coloneqq
        \{(r,z)\in\bA_1\times\bG_m:z^{2N}h_{g_{f'}}(r)=1\}.
    \]
    The projection
    \[
        \widetilde D_{f'}\longrightarrow D_{f'},\qquad
        (r,z)\longmapsto y_{f'}(r)
    \]
    is finite \'etale of degree $2N$, since $2N$ is prime to $\opr{char}(k)$.
    On $\widetilde D_{f'}$, write
    $(x,y)^{\mathsf t}=g_{f'}(r,1)^{\mathsf t}$ and set
    \[
        A(r,z)\coloneqq z^N\begin{pmatrix}x&y\\y&x\end{pmatrix}.
    \]
    The matrix $A(r,z)$ belongs to the fixed torus and
    \[
        fp_{f'}u_\infty(r)\dot sU
        =f\phi_{s,f}(A(r,z))\alpha_s^\vee(z)^{-N}U.
    \]
    Equivariance trivializes $(f\tilde y_{f'})^*\mathscr E_\xi$
    after pullback to $\widetilde D_{f'}$.
    The deck transformation $z\mapsto\zeta z$, for $\zeta^{2N}=1$,
    changes $A$ to $\zeta^N A$.
    Its character on the coefficient line is
    $\mu(\alpha_s^\vee(\zeta))$ by
    Lemma~\ref{lem:combi-parameter-rel-standard-obj}.
    Lemma~\ref{lem:rank-one-coefficient-restriction}\ref{itm:rank-one-geometric-coefficients} transfers this
    description to $\mathscr F_{f'}$.

    A zero or pole of $h_{g_{f'}}$ of order $m$ therefore has monodromy
    $\eta_{s,\xi}^m$, up to a common inversion.
    Direct evaluation gives
    \[
        h_1(r)=r^2-1,\qquad h_c(r)=-2r,\qquad
        h_{c\left(\begin{smallmatrix}0&-1\\1&0\end{smallmatrix}\right)}(r)=2r.
    \]
    Thus $h_1$ has two simple zeros and a double pole at $z_f=B$.
    Each boundary polynomial has one simple zero and one simple pole,
    giving the second row.
    A tame rank-one system on a punctured projective line is
    constant exactly when all local monodromies are trivial.
\end{proof}

\begin{lemma}
    \label{lem:real-cohomology-characters}
    The point $\dot sB\in D_f$ is fixed by
    $T_{f,N}$. Its stalk character and the boundary permutation
    sign are
    \[
        \rho_x(t)=\mu(s(t)),\qquad
        \delta(t)=\alpha_s(t^N),\qquad
        \rho_x\delta=s\times_f\mu,\qquad
        t\in T_{f,N}.
    \]
    In the parity case, $\rho_x=\mu$.
    For a boundary output in the parity case, the character on
    $H^0(D_{f'},\mathscr F_{f'})$ is the unique member of
    $c^s(\xi)$ on that orbit.
\end{lemma}
\begin{proof}
    The first assertions follow from
    Lemmas~\ref{lem:rank-one-coefficient-restriction}\ref{itm:rank-one-stabilizer-character},
    \ref{lem:rank-one-boundary-permutation-sign},
    \ref{lem:cross-action-well-defined}, and
    \ref{lem:reflection-on-parameter-characters}.

    For a boundary output in the parity case, $\mathscr F_{f'}$ is
    geometrically constant by Lemma~\ref{lem:real-root-monodromy}.
    Since $T_{f',N}^\circ$ acts transitively on $D_{f'}\cong\bG_m$,
    Lemma~\ref{lem:global-section-character-from-stabilizer} reduces
    the identification of the global-section character to a point stabilizer.
    The output chart identifies the point stabilizer with
    $\{t\in T_{f',N}:\alpha_s(t^N)=1\}$,
    the common subgroup of Lemma~\ref{lem:cayley-common-stabilizer}
    applied to the imaginary boundary orbit.
    Lemma~\ref{lem:rank-one-coefficient-restriction}\ref{itm:rank-one-stabilizer-character} gives its stalk character $\mu\circ s$.
    Real parity gives $\mu\circ s=\mu$ on this subgroup.
    The restriction condition in Definition~\ref{def:cayley-transform}
    and the uniqueness proved in Lemma~\ref{lem:cayley-transform-computation}
    identify the resulting boundary character.
\end{proof}

\begin{proposition}
    \label{prop:real-root-action}
    For a real root,
    \[
        \mathbf T_sm_\xi=
        \begin{cases}
            (u-1)m_\xi-m_{s\times\xi}
            +(u-1)\displaystyle\sum_{\zeta\in c^s(\xi)}m_\zeta,
                                   & \text{parity},                              \\
            -m_\xi,                & \text{non-parity and }s\cL_\xi\cong\cL_\xi, \\
            u^{1/2}m_{s\times\xi}, & s\cL_\xi\not\cong\cL_\xi.
        \end{cases}
    \]
\end{proposition}
\begin{proof}
    In the parity case, every coefficient system is constant by
    Lemma~\ref{lem:real-root-monodromy}.
    At $f'=f$, Lemma~\ref{lem:action-on-cohomology} gives
    \[
        H_c^1(D_f,\mathscr F_f)\cong\mu\oplus s\times_f\mu,
        \qquad H_c^2(D_f,\mathscr F_f)\cong\mu(-1).
    \]
    Their weights are zero and two, respectively.
    The contribution is $(u-1)m_\xi-m_{s\times\xi}$.
    Each boundary output has domain $\bG_m$ and contributes
    $u-1$ times its inverse Cayley parameter by
    Lemma~\ref{lem:real-cohomology-characters}.

    In the non-parity case, every two-puncture domain has zero
    cohomology. On $D_f$, only $H_c^1$ survives.
    Its character is $\rho_x\delta=s\times_f\mu$ by
    Lemma~\ref{lem:action-on-cohomology}.
    If $s\cL_\xi\cong\cL_\xi$, then
    $\eta_{s,\xi}=-1$ and $s\times\xi=\xi$.
    Exactly one local monodromy is trivial, so $H_c^1$ has weight
    zero and contributes $-m_\xi$.
    Otherwise $\eta_{s,\xi}^2\ne1$, and all three monodromies
    are nontrivial. The cohomology line has weight one, so its
    coefficient is $(-1)^1(-u^{1/2})^1=u^{1/2}$.
\end{proof}

\subsection{The complete action formula}
\label{subsec:complete-action-formula}
\begin{theorem}
    \label{prop:rank-one-action-formulas}
    Let \(s\in S\) and
    \(\xi=(f,\mu)\in\Gamma^{H^\circ}_\fo\).
    In the table, $\xi_s$ denotes the unique member, and
    $\xi_s^+,\xi_s^-$ the two members, of $c_s(\xi)$ in the
    non-compact imaginary rows and of $c^s(\xi)$ in the real parity rows.
    The labels $+$ and $-$ are arbitrary.
    \begin{enumerate}
        \item If \(s\cL_\xi\cong\cL_\xi\), then
              \[
                  \renewcommand{\arraystretch}{1.25}
                  \begin{array}{c|c|l}
                      \text{Root type} & \text{Case}           & \mathbf T_s m_\xi             \\ \hline
                      \text{complex}   & \text{ascent}         & m_{s\times\xi}                \\
                      \text{complex}   & \text{descent}        & (u-1)m_\xi+u m_{s\times\xi}   \\
                      \text{imaginary} & \text{compact}        & u m_\xi                       \\
                      \text{imaginary} & \text{non-compact I}  & m_{s\times\xi}+m_{\xi_s}      \\
                      \text{imaginary} & \text{non-compact II} & m_\xi+m_{\xi_s^+}+m_{\xi_s^-} \\
                      \text{real}      & \text{parity I}       &
                      (u-2)m_\xi+(u-1)(m_{\xi_s^+}+m_{\xi_s^-})                                \\
                      \text{real}      & \text{parity II}      &
                      (u-1)m_{\xi_s}+(u-1)m_\xi-m_{s\times\xi}                                 \\
                      \text{real}      & \text{non-parity}     & -m_\xi.
                  \end{array}
              \]
        \item If \(s\cL_\xi\not\cong\cL_\xi\), the pair is complex or real, and
              \[
                  \mathbf T_s m_\xi=
                  \begin{cases}
                      m_{s\times\xi},        & \text{complex ascent},  \\
                      u m_{s\times\xi},      & \text{complex descent}, \\
                      u^{1/2}m_{s\times\xi}, & \text{real}.
                  \end{cases}
              \]
    \end{enumerate}
\end{theorem}
\begin{proof}
    The complex, imaginary, and real rows are proved in
    Propositions~\ref{prop:complex-root-action},
    \ref{prop:imaginary-root-action}, and~\ref{prop:real-root-action},
    respectively. Lemma~\ref{lem:cayley-transform-computation}
    gives the displayed cardinalities and cross-action identities.
\end{proof}

\section{Specialization and the leading-term formula}
\label{sec:specialization-leading-term}
We specialize the operator formulas obtained from the rank-one
models. The cross-action term preserves the anti-invariant torus
rank, while every surviving additional term increases it.
\subsection{Rank and the specialized action}
\label{subsec:specialized-action}
Retain $M^{H^\circ}_{\fo}$ and the specialization
$\epsilon(u^{1/2})=-1$ of
\S\ref{subsubsec:the-generic-module-and-frobenius}.
The following degree measures the anti-invariant part of the torus.

\begin{definition}
    \label{def:theta-degree}
    For $\xi=(f,\mu)\in\Gamma^{H^\circ}_{\fo}$, put
    \[
        T^{-\theta_f}\coloneqq\{t\in T\mid\theta_f(t)=t^{-1}\},
        \qquad r_\theta(\xi)\coloneqq\dim (T^{-\theta_f})^\circ.
    \]
\end{definition}

\begin{lemma}
    \label{lem:theta-degree-and-rank-one-transforms}
    For every \(w\in W\), \(r_\theta(w\times\xi)=r_\theta(\xi)\).
    If \((s,\xi)\) is non-compact imaginary and
    \(\zeta\in c_s(\xi)\), then
    \[
        r_\theta(\zeta)=r_\theta(\xi)+1.
    \]
\end{lemma}
\begin{proof}
    Write $\xi=(f,\mu)$.
    Lemma~\ref{lem:cross-action-well-defined} gives
    $T^{-\theta_{w\times f}}=w(T^{-\theta_f})$, proving the first claim.

    Suppose $(s,\xi)$ is non-compact imaginary and
    $\zeta=(f',\mu')\in c_s(\xi)$.
    Definition~\ref{def:cayley-transform} gives $f'=f_s$, and
    Lemma~\ref{lem:rank-one-cayley-geometry} gives
    $\theta_{f'}|_T=s\circ\theta_f|_T$.

    The decomposition
    $X_*(T)_{\bQ}=\ker(\alpha_s)\oplus\bQ\alpha_s^\vee$
    is $\theta_f$-stable, and $\theta_f$ fixes $\alpha_s^\vee$.
    The reflection $s$ fixes $\ker(\alpha_s)$ pointwise and negates
    $\alpha_s^\vee$. Thus
    \[
        X_*(T)_{\bQ}^{\theta_{f'}=-1}
        =X_*(T)_{\bQ}^{\theta_f=-1}
        \oplus\bQ\alpha_s^\vee.
    \]
    Since $r_\theta(\xi)=\dim_{\bQ}X_*(T)_{\bQ}^{\theta_f=-1}$,
    this proves $r_\theta(\zeta)=r_\theta(\xi)+1$.
\end{proof}

\begin{lemma}
    \label{lem:specialized-rank-one-leading-term}
    Let $s\in S$ and $\xi=(f,\mu)\in\Gamma^{H^\circ}_{\fo}$.
    With $\epsilon(u^{1/2})=-1$, one has
    \[
        \epsilon(\mathbf T_s)m_\xi-c_{s,\xi}m_{s\times\xi}
        \in\opr{span}_{\overline{\bQ}_\ell}
        \{m_\zeta\mid\zeta\in\Gamma^{H^\circ}_{\fo},\
        r_\theta(\zeta)>r_\theta(\xi)\},
    \]
    where
    \[
        c_{s,\xi}\coloneqq
        \begin{cases}
            -1, & \theta_f(\alpha_s)=-\alpha_s, \\
            1,  & \text{otherwise}.
        \end{cases}
    \]
\end{lemma}

\begin{proof}
    Specialize Theorem~\ref{prop:rank-one-action-formulas}
    at $u^{1/2}=-1$. The cross-action coefficient is $-1$
    for real roots and $1$ otherwise, including the cases
    $s\cL_\xi\not\cong\cL_\xi$.
    Every remaining term is indexed by $\zeta\in c_s(\xi)$
    for a non-compact imaginary input.
    Lemma~\ref{lem:theta-degree-and-rank-one-transforms} gives
    $r_\theta(\zeta)=r_\theta(\xi)+1$ and
    $r_\theta(s\times\xi)=r_\theta(\xi)$.
\end{proof}

\begin{corollary}
    \label{cor:rank-filtration}
    For every integer $r$, set
    \[
        \mathrm{Fil}^{\ge r}M^{H^\circ}_{\fo}
        \coloneqq\opr{span}_{\overline{\bQ}_\ell}
        \{m_\xi\mid\xi\in\Gamma^{H^\circ}_{\fo},\ r_\theta(\xi)\ge r\}.
    \]
    These subspaces form a $W$-stable decreasing filtration, with
    $\mathrm{Fil}^{\ge0}M^{H^\circ}_{\fo}=M^{H^\circ}_\fo$ and
    $\mathrm{Fil}^{\ge\rank G+1}M^{H^\circ}_{\fo}=0$.
    On $\mathrm{Fil}^{\ge r}M^{H^\circ}_{\fo}/
    \mathrm{Fil}^{\ge r+1}M^{H^\circ}_{\fo}$,
    the basis classes with $r_\theta(\xi)=r$ satisfy
    \[
        s\cdot\overline m_\xi=c_{s,\xi}\overline m_{s\times\xi},
        \qquad s\in S,
    \]
    where $c_{s,\xi}$ is the sign in
    Lemma~\ref{lem:specialized-rank-one-leading-term}.
\end{corollary}


\begin{thebibliography}{BBD82}

\bibitem[BBD82]{BBD82}
A.~A. Be\u{\i}linson, J.~Bernstein, and P.~Deligne.
\newblock Faisceaux pervers.
\newblock In {\em Analysis and topology on singular spaces, {I} ({L}uminy,
  1981)}, volume 100 of {\em Ast\'{e}risque}, pages 5--171. Soc. Math. France,
  Paris, 1982.

\bibitem[Dag26]{dagan-distinction}
Shahar Dagan.
\newblock {Lusztig} families of distinguished representations for symmetric
  spaces.
\newblock Preprint, 2026.

\bibitem[DL76]{deligneLusztig1976}
Pierre Deligne and George Lusztig.
\newblock Representations of reductive groups over finite fields.
\newblock {\em Ann. of Math. (2)}, 103(1):103--161, 1976.

\bibitem[Hel94]{helminck1994symmetric}
A.~G. Helminck.
\newblock Symmetric {$k$}-varieties.
\newblock In {\em Algebraic groups and their generalizations: Classical
  methods. {I}}, volume~56 of {\em Proceedings of Symposia in Pure
  Mathematics}, pages 233--279. American Mathematical Society, Providence, RI,
  1994.

\bibitem[Hen01]{henderson2003character}
Anthony Henderson.
\newblock {\em Character sheaves on symmetric spaces}.
\newblock PhD thesis, Massachusetts Institute of Technology, 2001.

\bibitem[HW93]{helminck1993rationality}
A.~G. Helminck and S.~P. Wang.
\newblock On rationality properties of involutions of reductive groups.
\newblock {\em Adv. Math.}, 99(1):26--96, 1993.

\bibitem[Kat13]{katz2013squarefree}
Nicholas~M. Katz.
\newblock On a question of {Keating} and {Rudnick} about primitive {Dirichlet}
  characters with squarefree conductor.
\newblock {\em Int. Math. Res. Not. IMRN}, (14):3221--3249, 2013.

\bibitem[Lus19]{lusztig2017conjugacy}
George Lusztig.
\newblock Conjugacy classes in reductive groups and two-sided cells.
\newblock {\em Bull. Inst. Math. Acad. Sin. (N.S.)}, 14(3):265--293, 2019.

\bibitem[LV83]{lusztig1983singularities}
George Lusztig and David~A. Vogan, Jr.
\newblock Singularities of closures of {K}-orbits on flag manifolds.
\newblock {\em Invent. Math.}, 71(2):365--379, 1983.

\bibitem[LV14]{lusztig-vogan2014quasisplit}
George Lusztig and David~A. Vogan, Jr.
\newblock Quasisplit {Hecke} algebras and symmetric spaces.
\newblock {\em Duke Math. J.}, 163(5):983--1034, 2014.
\newblock References to sections follow arXiv:1206.0634v4.

\bibitem[LY20]{lusztig2020endoscopy}
George Lusztig and Zhiwei Yun.
\newblock Endoscopy for {Hecke} categories, character sheaves and
  representations.
\newblock {\em Forum Math. Pi}, 8:e12, 2020.

\bibitem[Mil13]{milne2013etale}
James~S. Milne.
\newblock Lectures on {\'e}tale cohomology, 2013.
\newblock Version 2.21.

\end{thebibliography}
\end{document}